\documentclass[12pt]{article}
\usepackage{amsmath}
\usepackage{amssymb}
\usepackage{amsfonts}
\usepackage{graphicx}
\usepackage[all]{xy}

\begin{document}

\title{\textbf{\normalsize {\ SEMI-CUBIC HYPONORMALITY OF WEIGHTED SHIFTS WITH STAMPFLI RECURSIVE TAIL}}}
\author{{\normalsize Seunghwan Baek, George R. Exner, Il Bong Jung and
Chunji Li\thanks{%
The research of the third author was supported by Basic Science Research
Program through the National Research Foundation of Korea(NRF) funded by the
Ministry of Science, ICT and Future Planning (KRF-2015R1A2A2A01006072).}}}
\date{}
\maketitle

\begin{abstract}
Let $\alpha :\sqrt{x_{m}},\cdots ,\sqrt{x_{1}},(\sqrt{u},\sqrt{v},\sqrt{w}%
)^{\wedge}$ be a backward $m$-step extension of a recursive weight sequence and let $%
W_{\alpha }$ be the weighted shift associated with $\alpha $. In this paper we
characterize the semi-cubic hyponormality of $W_{\alpha }$ having the positive
determinant coefficient property and discuss some related examples.
\end{abstract}

\textbf{1. Introduction and Notation. }Let $\mathcal{H}$ be a separable
infinite dimensional complex Hilbert space and let $\mathcal{L}(\mathcal{H})$
be the algebra of all bounded linear operators on $\mathcal{H}$. For $X$, $%
Y\in \mathcal{L}(\mathcal{H})$, we set $[X,Y]:=XY-YX$. An operator $T$ in
$\mathcal{L(H)}$ is \textit{subnormal} if it is (unitarily equivalent to)
the restriction of a normal operator to an invariant subspace, and \textit{%
hyponormal} if $[T^{\ast },T]\geq 0$. There are several classes of operators in $%
\mathcal{L(H)}$ between subnormal and hyponormal which have been discussed for more than two decades (cf. \cite{BEJ}, \cite{Cu1}, \cite{Cu2}, \cite{CuF1}, \cite{CuF2}, \cite{EJPa}, \cite{JuP2}), including,  for
example, (strongly) $n$-hyponormal and weakly $n$-hyponormal operators. Some of these
definitions were derived from the theory of multivariable operators. A $k$%
-tuple $\mathbf{T}=(T_{1},...,T_{k})$ of operators on $\mathcal{H}$ is
called \textit{hyponormal} if the operator matrix $([T_{j}^{\ast
},T_{i}])_{i,j=1}^{k}$ is positive on the direct sum of $\mathcal{H}\oplus
\cdots \oplus \mathcal{H}$ with $k$ copies. An operator $T$ is said to be (%
\textit{strongly}) $k$-\textit{hyponormal} if $(I,T,...,T^{k})$ is
hyponormal for $k\in \mathbb{N}$, where $\mathbb{N}$ is the set of positive
integers. It is known that $T\in \mathcal{L(H)}$ is subnormal if and only if
$T$ is $k$-hyponormal for all $k\in \mathbb{N}$ (\cite{Br}, \cite{Hal}, \cite{Con}). An
operator $T$ in $\mathcal{L}(\mathcal{H})$ is \textit{weakly} $k$\textit{%
-hyponormal} if $p(T)$ is hyponormal for every polynomial $p$ of degree $k$
or less (\cite{CuF1}, \cite{CuP1}, \cite{CuP2});  it is \textit{polynomially hyponormal} if $p(T)$ is hyponormal for each polynomial $p$.   In particular, weakly $2$-hyponormal (or weakly $3$%
-hyponormal) is referred to as \textit{quadratically hyponormal} (or \textit{%
cubically hyponormal}, resp.).  For $k\in \mathbb{N}$, an
operator $T\in \mathcal{L}(\mathcal{H})$ is called \textit{semi-weakly} $k$%
\textit{-hyponormal} if $T+sT^{k}$ is hyponormal for all $s\in \mathbb{C}$ (%
\cite{DEJL}).  It is obvious that a weakly $k$-hyponormal
operator is semi-weakly $k$-hyponormal. Note also that $k$-hyponormality implies
weak $k$-hyponormality for $k\in \mathbb{N}$. The following implications
provide a bridge between subnormal and hyponormal operators:

\begin{center}
subnormal $\Rightarrow $ polynomially hyponormal $\Rightarrow \cdots
\Rightarrow $ weakly $3$-hyponormal\\
$\Rightarrow $ weakly $2$-hyponormal $%
\Rightarrow $ hyponormal.
\end{center}

\noindent In \cite{CuP1} and \cite{CuP2}, Curto-Putinar proved the existence
of polynomially hyponormal weighted shifts which are not $2$-hyponormal. But
no concrete example of such a weighted shift has been found for twenty
five years. Quadratic hyponormality and cubic hyponormality are
related closely to attempts to construct such examples. In a different line of study, J. Stampfli (%
\cite{St}) constructed the subnormal completion induced by $(\sqrt{u},\sqrt{v%
},\sqrt{w})^{\wedge}$ with $0<u<v<w$. In \cite{CuF2}, Curto-Fialkow
characterized the positive quadratic hyponormality of a weighted shift $%
W_{\alpha }$ associated to a backward 1-step extension $\alpha :\sqrt{x},(%
\sqrt{u},\sqrt{v},\sqrt{w})^{\wedge}.$ In \cite{JuP1}, Jung-Park improved
their results about such weighted shifts; namely, it was proved that $%
W_{\alpha }$ is quadratically hyponormal if and only if $W_{\alpha }$ is
positively quadratically hyponormal (see that paper for the relevant definition). In \cite{EJLP}, one characterized the
positive quadratic hyponormality of a weighted shift $W_{\alpha }$
associated to a backward $m$-step extension $\alpha :\sqrt{x_{m}},\cdots ,%
\sqrt{x_{1}},(\sqrt{u},\sqrt{v},\sqrt{w})^{\wedge}.$ In \cite{PY}, Poon-Yoon
obtained a weighted shift $W_{\alpha }$ associated to a backward $2 $-step
extension $\alpha :\sqrt{x_{2}},\sqrt{x_{1}},(\sqrt{u},\sqrt{v},\sqrt{w}%
)^{\wedge},$ which is quadratically hyponormal but not positively
quadratically hyponormal. In yet another direction of study, J. Stampfli (\cite{St}) proved that a
subnormal weighted shift with two equal weights $\alpha _{n}=\alpha _{n+1}$
for some nonnegative $n$ has the flatness property, i.e., $\alpha
_{1}=\alpha _{2}=\cdots $. Stampfli's result has been used to attempt the
construction of nonsubnormal polynomially hyponormal weighted shifts (cf.
\cite{Ch}, \cite{Cu1}, \cite{DEJL}, \cite{JuP1}). In \cite{Ch}, it was proved
that if a weighted shift $W_{\alpha }$ is polynomially hyponormal with first
two weights equal, then $W_{\alpha }$ has flatness. The notion of semi-weak
hyponormality of a weighted shift $W_{\alpha }$ comes from the study of
polynomial hyponormality via the positive determinant coefficients property (cf.
\cite{DEJL}, or see below). In this paper we discuss semi-cubic
hyponormality of a weighted shift $W_{\alpha }$ associated to a backward $m$%
-step extension $\alpha :\sqrt{x_{m}},\cdots ,\sqrt{x_{1}},(\sqrt{u},\sqrt{v}%
,\sqrt{w})^{\wedge}$ via the positive determinant coefficients property.  Here we consider a recursive tail, which may be contrasted with \cite{DEJL} which considered related questions but with the tail arising from the Bergman shift (having weights $\sqrt{\frac{1}{2}}, \sqrt{\frac{2}{3}}, \sqrt{\frac{3}{4}}, \ldots$, which sequence is known to be non-recursive).

The paper consists of the following. In Section 2 we recall some terminology
concerning semi-cubically hyponormal weighted shifts. In Section 3 we
characterize the positive determinant coefficients property for
semi-cubically hyponormal weighted shifts. In Section 4 we characterize the
semi-cubic hyponormality of such weighted shifts $W_{\alpha }$ with the positive
determinant coefficients property, where $\alpha :\sqrt{x_{m}},\cdots ,\sqrt{%
x_{1}},(\sqrt{u},\sqrt{v},\sqrt{w})^{\wedge}.$ In Section 5, we discuss some
related examples;  we omit many of the details of these computations, but they may be found in an extended version of this paper (particularly in the Appendices) posted on ar{$\chi$}iv.  Many of these computations, and the included graphs, were obtained with the aid of \textit{Mathematica}$^{\copyright}$ (\cite{Wol}).

Throughout this paper, $\mathbb{C},\mathbb{R}_{+},\mathbb{N}$, and $\mathbb{N}_{0}$ are the
sets of complex, nonnegative real numbers, positive integers, and nonnegative
integers, respectively.

\bigskip

\textbf{2. Preliminaries.} We recall some standard terminology for
semi-cubically hyponormal weighted shifts (cf. \cite{DEJL}). Let $\ell ^{2}(%
\mathbb{N}_{0})$ be the space of square summable sequences in $\mathbb{C}%
_{+} $ and let $\{e_{i}\}_{i=0}^{\infty }$ be the standard orthonormal basis of $\ell
^{2}(\mathbb{N}_{0})$. For a weight sequence $\alpha =\{\alpha
_{i}\}_{i=0}^{\infty }$ in $\mathbb{R}_{+}$, the associated weighted shift $%
W_{\alpha }$ acting on $\ell ^{2}(\mathbb{N}_{0})$ is defined by $W_\alpha e_j = \alpha_j e_{j+1}$ and extended by linearity.  The shift $W_\alpha$ is \textit{semi-cubically
hyponormal} if
\begin{equation}
D(s):=[(W_{\alpha }+sW_{\alpha }^{3})^{\ast },W_{\alpha }+sW_{\alpha
}^{3}]\geq 0,\text{ \ \ }s\in \mathbb{C}.  \tag{2.1}
\end{equation}%
\medskip In fact, the condition (2.1) is equivalent to the (\textit{a priori} weaker) condition that $W_{\alpha
}+tW_{\alpha }^{3}$\ is hyponormal for all $t\geq 0$ (cf. (\cite[Prop. 1]%
{CuL}, \cite[Prop. 2.1]{BEJL}). Hence $W_{\alpha }$ is semi-cubically
hyponormal if and only if $D(s)\geq 0$ for all $s\in \mathbb{R}_{+}$.
Observe that
\begin{equation*}
D(s)=\left(
\begin{array}{ccccc}
q_{0} & 0 & r_{0} & 0 & \cdots \\
0 & q_{1} & 0 & r_{1} & \ddots \\
r_{0} & 0 & q_{2} & \ddots & \ddots \\
0 & r_{1} & \ddots & \ddots & \ddots \\
\vdots & \ddots & \ddots & \ddots & \ddots%
\end{array}%
\right) ,~\ \ s\in \mathbb{R}_{+},
\end{equation*}%
where for $k\in \mathbb{N}_{0},$
\begin{align*}
q_{k}& :=u_{k}+v_{k}s^{2},\ \ r_{k}:=\sqrt{w_{k}}s,\ \ u_{k}:=\alpha
_{k}^{2}-\alpha _{k-1}^{2}, \\
v_{k}& :=\alpha _{k}^{2}\alpha _{k+1}^{2}\alpha _{k+2}^{2}-\alpha
_{k-3}^{2}\alpha _{k-2}^{2}\alpha _{k-1}^{2},\ w_{k}:=\alpha _{k}^{2}\alpha
_{k+1}^{2}(\alpha _{k+2}^{2}-\alpha _{k-1}^{2})^{2}
\end{align*}%
with $\alpha _{-3}=\alpha _{-2}=\alpha _{-1}=0$. Considering two submatrices
\begin{equation*}
D^{[1]}(s)=\left(
\begin{array}{lllll}
q_{0} & r_{0} & 0 &  &  \\
r_{0} & q_{2} & r_{2} & 0 &  \\
0 & r_{2} & q_{4} & r_{4} & \ddots \\
& 0 & r_{4} & \ddots & \ddots \\
&  & \ddots & \ddots & \ddots%
\end{array}%
\right) \text{ and }D^{[2]}(s)=\left(
\begin{array}{lllll}
q_{1} & r_{1} & 0 &  &  \\
r_{1} & q_{3} & r_{3} & 0 &  \\
0 & r_{3} & q_{5} & r_{5} & \ddots \\
& 0 & r_{5} & \ddots & \ddots \\
&  & \ddots & \ddots & \ddots%
\end{array}%
\right) ,
\end{equation*}%
we have a decomposition $D(s)=D^{[1]}(s)\oplus D^{[2]}(s),$ $s\in \mathbb{R}%
_{+}.$ Define
\begin{equation}
D_{n}^{[1]}(t)=\left(
\begin{array}{lllll}
q_{0} & r_{0} & 0 &  &  \\
r_{0} & q_{2} & r_{2} & 0 &  \\
0 & r_{2} & q_{4} & r_{4} & \ddots \\
& 0 & r_{4} & \ddots & \ddots \\
&  & \ddots & \ddots & q_{2n}%
\end{array}%
\right) ,\ \ D_{n}^{[2]}(t)=\left(
\begin{array}{lllll}
q_{1} & r_{1} & 0 &  &  \\
r_{1} & q_{3} & r_{3} & 0 &  \\
0 & r_{3} & q_{5} & r_{5} & \ddots \\
& 0 & r_{5} & \ddots & \ddots \\
&  & \ddots & \ddots & q_{2n+1}%
\end{array}%
\right) ,  \tag{2.2}
\end{equation}%
where $t=s^{2}$. Then $W_{\alpha }$ is semi-cubically hyponormal if and only
if $D_{n}^{[j]}(t)\geq 0$ for all $n\geq 0,$ $j=1,2$.

To consider the two matrices $D_{n}^{[j]}(t)$ in (2.2), we consider the matrix
with the form below%
\begin{equation}  \label{def:genmatrix}
M_{n}\left( t\right) =\left(
\begin{array}{cccccc}
\check{q}_{0} & \check{r}_{0} & 0 &  &  &  \\
\check{r}_{0} & \check{q}_{1} & \check{r}_{1} & 0 &  &  \\
0 & \check{r}_{1} & \check{q}_{2} & \check{r}_{2} & \ddots &  \\
& 0 & \check{r}_{2} & \ddots & \ddots & 0 \\
&  & \ddots & \ddots & \check{q}_{n-1} & \check{r}_{n-1} \\
&  &  & 0 & \check{r}_{n-1} & \check{q}_{n}%
\end{array}%
\right) ,  \tag{2.3}
\end{equation}%
where $\check{q}_{k}:=\check{u}_{k}+\check{v}_{k}t,r_{k}:=\sqrt{\check{w}%
_{k}t}$ $\left( k\geq 0\right) ,$ and $\check{u}_{k}\geq 0,\check{v}_{k}\geq
0,\check{w}_{k}\geq 0,t\geq 0.$ Set
\begin{equation}
d_{n}\left( t\right) :=\det M_{n}(t)=\sum_{i=0}^{n+1}c\left( n,i\right)
t^{i};  \tag{2.4}
\end{equation}%
it follows from \cite{CuF2} that
\begin{align*}
c(0,0)& =\check{u}_{0},\ \ c(0,1)=\check{v}_{0}, \\
c(1,0)& =\check{u}_{0}\check{u}_{1},\ c(1,1)=\check{u}_{1}\check{v}_{0}+%
\check{u}_{0}\check{v}_{1}-\check{w}_{0},\ c(1,2)=\check{v}_{1}\check{v}_{0},
\\
c(n+2,i)& =\check{u}_{n+2}c(n+1,i)+\check{v}_{n+2}c(n+1,i-1)-\check{w}%
_{n+1}c(n,i-1), \\
c(n,n+1)& =\check{v}_{0}\check{v}_{1}\cdots \check{v}_{n}, \\
c(n,0)& =\check{u}_{0}\check{u}_{1}\cdots \check{u}_{n},\ \ n\geq 0,0\leq
i\leq n,
\end{align*}%
with $c(-n,-i):=0$ for all $n,i\in \mathbb{N}$.  We say that the matrix in (\ref{def:genmatrix}) has \textit{positive determinant coefficients} if each $c(n,i)$ is non-negative and $c(n, n+1) > 0$ for each $n$.

This definition is in anticipation of the use of the Nested Determinant Test.  Recall that for a real, symmetric matrix $M$, if the determinants of the principal submatrices of $M$ are positive and det $M \geq0$, then $M$ is positive semi-definite.  The assumption that $c(n, n+1) > 0$ should be regarded as a convenient way to achieve the that the submatrices have strictly positive determinants, and in our applications will involve no essential loss of generality.  See the discussion in the Remark at the beginning of Section 4 for further details.

\bigskip

\textbf{3. A characterization of positive determinant coefficients.}
In what follows in this section we will make the additional assumption that $\check{v}_k > 0$ for all $k$.  In light of the formula for $c(n,n+1)$ above this assumption achieves that requirement of the definition of positive determinant coefficients, and we will omit further reference to it.

We
first consider a new matrix related to $M_n(t)$ as in (\ref{def:genmatrix}) as follows:
\begin{equation*}
\widehat{M}_{n}(t):=M_{n}(t)-\left(
\begin{array}{ccccc}
0 & 0 & \cdots & 0 & 0 \\
0 & 0 & \cdots & 0 & 0 \\
\vdots & \vdots & \ddots & \vdots & \vdots \\
0 & 0 & \cdots & 0 & 0 \\
0 & 0 & \cdots & 0 & \check{u}_{n}%
\end{array}%
\right) .
\end{equation*}%
Then we get obviously that det$\widehat{M}_{n}(t)=$det$M_{n}(t)-\check{u}%
_{n} $det$M_{n-1}(t)$ with det$\widehat{M}_{n}(t)$ having a Maclaurin expansion
\begin{equation*}
\text{det}\widehat{M}_{n}(t)=\sum_{i=0}^{n+1}\widehat{c}(n,i)t^{i}.
\end{equation*}%
Observe that $\widehat{c}(n,n+1)=c(n,n+1)$ and also, for $i=0$, we
have
\begin{equation*}
\widehat{c}(n,0)=c(n,0)-\check{u}_{n}c(n-1,0)=0.
\end{equation*}%
Further, the recursion in the $c(n,i)$ induces that, for $1\leq i\leq n$ and $n\geq 2$,
\begin{align*}
\widehat{c}(n,i)& =c(n,i)-\check{u}_{n}c(n-1,i) \\
& =\check{v}_{n}c(n-1,i-1)-\check{w}_{n-1}c(n-2,i-1).
\end{align*}%
As will be true in our applications, we assume that for some positive integer $m_0$%
\begin{equation}
\check{u}_{n}\check{v}_{n+1}=\check{w}_{n},\text{ }n\geq m_{0}.  \tag{3.1}
\end{equation}%

In what follows we imitate the argument in \cite{EJLP}.
Suppose $n\geq m_{0}+1$ for the time being
and set $\check{\eta}_{n}=\frac{\check{v}_{n}}{\check{u}_{n}}$. Then
we have a coefficient formula (cf. \cite{EJLP}) (note that what we record next does not include all possible pairs $(n,i)$):%
\begin{equation}
c(n,i)=\left\{
\begin{array}{ll}
\check{u}_{n}c(n-1,i), & 0\leq i\leq n-m_{0}, \\
\check{u}_{n}c(n-1,i)+\check{v}_{n}\cdots \check{v}_{m_{0}+1}\widehat{c}%
(m_{0},i-n+m_{0}), & n-m_{0}+1\leq i\leq m_{0}, \\
\check{u}_{m_{0}+1}(\check{v}_{0}\cdots \check{v}_{m_{0}}+\check{\eta}_{n}%
\widehat{c}(m_{0},m_{0})), & i=n=m_{0}+1, \\
\check{u}_{n}\check{v}_{n-1}\cdots \check{v}_{m_{0}+1}(\check{v}_{0}\cdots
\check{v}_{m_{0}}+\check{\eta}_{n}\widehat{c}(m_{0},m_{0})), & i=n>m_{0}+1,
\\
\check{v}_{n}\cdots \check{v}_{0}, & i=n+1.%
\end{array}%
\right.  \tag{3.2}
\end{equation}%
Note that the sign of $c(n,i)$ $(0\leq i\leq m_{0}$, $i=n,n+1)$ can be
checked using this coefficient formula.

For the coefficient pairs not listed above (and which correspond to the numbered diagonals in Figure 1), it turns out to be convenient to write them as $c(n,n-i)$, because the relevant conditions on $n$ and $i$ are that both  $m_{0}+1\leq n-i\leq n$ and $n-m_{0}<n-i\leq n$.  Equivalently, these are $m_{0}+1+i\leq n$ and $0\leq i\leq m_{0}-1$, and these cover a variety of cases depending on the relative sizes of $n$ and $m_0$.  Using the recurrence, the $c(n,n-i)$ for $m_{0}+1+i\leq n$ and $0\leq i\leq m_{0}-1,$ can be represented
by
\begin{align*}
c(n,n-i)& =\check{u}_{n}\cdots \check{u}_{n-i}\check{v}_{n-i-1}\cdots \check{%
v}_{m_{0}+1} \\
& \times \left( \check{v}_{m_{0}}\cdots \check{v}_{0}+\sum_{k=0}^{i}\left(
\prod_{j=0}^{k}\check{\eta}_{n-i+j}\right) \widehat{c}(m_{0},m_{0}-k)\right)
.
\end{align*}%
To check the sign of these $c(n,n-i)$, it is convenient to consider
\begin{align*}
\check{G}(n,i)& \equiv \check{G}(n,i,\check{\eta}_{n-i},\cdots ,\check{\eta}_{n}) \\
& :=\check{v}_{m_{0}}\cdots \check{v}_{0}+\sum_{k=0}^{i}\left(
\prod_{j=0}^{k}\check{\eta}_{n-i+j}\right) \widehat{c}(m_{0},m_{0}-k),
\end{align*}%
and define
\begin{align*}
\Delta \check{G}(n,i):=& \check{G}(n+1,i)-\check{G}(n,i) \\
=& \sum_{k=0}^{i}\left( \prod_{j=0}^{k}\check{\eta}_{n-i+j+1}-%
\prod_{j=-1}^{k-1}\check{\eta}_{n-i+j-1}\right) \widehat{c}(m_{0},m_{0}-k),
\end{align*}%
(again, all this for $m_{0}+1+i\leq n$ and $0\leq i\leq m_{0}-1$). For brevity, we set
\begin{equation}
Q(\ell,k):=\frac{\prod_{j=0}^{k}\check{\eta}_{\ell +j+1}-\prod_{j=-1}^{k-1}\check{%
\eta}_{\ell+j+1}}{\check{\eta}_{\ell+1}-\check{\eta}_{\ell}}.  \tag{3.3}
\end{equation}%
Then for $m_{0}+1+i\leq n$ and $0\leq i\leq m_{0}-1$, we have that%
\begin{equation*}
\Delta \check{G}(n,i)=(\check{\eta}_{n-i+1}-\check{\eta}_{n-i})\sum_{k=0}^{i}Q(n-i,k)%
\widehat{c}(m_{0},m_{0}-k).
\end{equation*}%
And we assume from now on that
\begin{equation*}
\check{\eta}_{n+1}>\check{\eta}_{n},\ n\geq m_{0}+1.
\end{equation*}%
Suppose now that for some $i$ with $0\leq i\leq m_{0}-1,$ $\check{G}(n,i)$ is
increasing in $n$ for $n\geq m_{0}+1+i.$ To check the positivity of all $%
c(n,n-i)$ for $n\geq m_{0}+1+i$, it is enough to check the sign of $%
\check{G}(m_{0}+1+i,i)$. If, on the other hand, for some $i$, $0\leq i\leq m_{0}-1$,
$\check{G}(n,i)$ is decreasing in $n$, $(m_{0}+1+i\leq n)$, then to check the
positivity of all $c(n,n-i)$ for $n\geq m_{0}+1+i$, it is enough to see if $%
\lim_{n\rightarrow \infty }\check{G}(n,i)\geq 0$ (in our applications this limit will exist and be finite). If, in particular, the sequence $%
\check{\eta}_{n}$ (assumed increasing) converges to some $\check{\eta} <\infty $, we
may simply examine the sign of
\begin{equation*}
\lim_{n\rightarrow \infty }\check{G}(n,i)=\check{v}_{m_{0}}\cdots \check{v}%
_{0}+\sum_{k=0}^{i}\check{\eta} ^{k+1}\widehat{c}(m_{0},m_{0}-k).
\end{equation*}%
Hence we arrive the following proposition.

\medskip

\textbf{Proposition 3.1.} \textit{With notation as above, fix} $i$ \textit{%
such that} $0\leq i\leq m_{0}-1$. \textit{Suppose that }$\check{G}(n,i)$ \textit{is monotone increasing or decreasing for }$n\geq m_0 + i +1$. \textit{Then the following conditions are
equivalent:}

(i) \textit{the coefficients} $c(n,n-i)$ \textit{satisfy} $c(n,n-i)\geq 0$ (\textit{all} $n\geq
m_{0}+i+1 $),

(ii)\textit{\ it holds that either}

\ \ \ \ (ii-a) $\check{G}(n,i)$ \textit{is increasing in} $n$, $n\geq m_{0}+i+1$
\textit{and} $\check{G}(m_{0}+1+i,i)\geq 0$,

\textit{or}

\ \ \ \ (ii-b) $\check{G}(n,i)$ \textit{is decreasing in} $n$, $n\geq m_{0}+i+1$
\textit{and} $\lim_{n\rightarrow \infty }\check{G}(n,i)\geq 0$.

\medskip

Note that, in particular, if $\lim_{n\rightarrow \infty }\check{\eta}%
_{n}=\check{\eta} <\infty $, then the second condition in (ii-b) may be represented
by
\begin{equation*}
\lim_{n\rightarrow \infty }\check{G}(n,i)=\check{v}_{m_{0}}\cdots \check{v}%
_{0}+\sum_{k=0}^{i}\check{\eta} ^{k+1}\widehat{c}(m_{0},m_{0}-k) \geq 0.
\end{equation*}

From the first term in (3.2), we get that for each $p\geq 1$, if we have $%
c(m_{0}-1+p,p)\geq 0$ then $c(m_{0}-1+p+j,p)\geq 0$ for all $j\geq 0$. It
follows from (2.4) that for all $n$, $c(n,n+1)\geq 0$ and $c(n,0)\geq 0$.
Thus we have the following result for positive determinant coefficients for
our matrices.

\medskip

\textbf{Proposition 3.2.} \textit{Under the above notation, assume that }$%
\check{u}_{n}\check{v}_{n+1}=\check{w}_{n}$ $(n\geq m_{0})$\textit{\ and }$%
\check{\eta}_{n+1}>\check{\eta}_{n}$ $(n\geq m_{0}+1)$.\textit{\ Then the
matrices }$M_{n}(t)$\textit{\ have positive determinant coefficients if and
only if the following assertions hold}

(i) \textit{for each }$i$\textit{, }$0\leq i\leq m_{0}-1$\textit{, one of
the conditions of Proposition 3.1 holds,}

(ii)\textit{\ for each }$n$\textit{, }$1\leq n\leq m_{0}$\textit{\ and each }%
$j$\textit{, }$1\leq j\leq n$\textit{, }$c(n,j)\geq 0$\textit{,}

(iii)\textit{\ for each }$n$\textit{, }$m_{0}+1\leq n\leq 2m_{0}-1$\textit{\
and each }$j$\textit{, }$n+1-m_{0}\leq j\leq m_{0}$\textit{, }$c(n,j)\geq 0$%
\textit{.}

\bigskip

Figure 1 gives an indication of the various groups of
coefficients $c(n,i)$.

\bigskip

\begin{figure}[ht]
\centerline{ \xy
(0,0);(0,88) **@{-}, (0,88);(88,88)**@{-},
(0,3);(1,3) **@{-}, (0,8);(1,8) **@{-},
(0,13);(1,13) **@{-}, (0,18);(1,18)**@{-},
(0,23);(1,23) **@{-}, (0,28);(1,28) **@{-},
(0,33);(1,33) **@{-}, (0,38);(1,38) **@{-},
(0,43);(1,43) **@{-}, (0,48);(1,48)**@{-},
(0,53);(1,53) **@{-}, (0,58);(1,58) **@{-},
(0,63);(1,63) **@{-}, (0,68);(1,68) **@{-},
(0,73);(1,73) **@{-}, (0,78);(1,78)**@{-},
(0,83);(1,83) **@{-}, (0,88);(1,88)**@{-},
(5,88);(5,87) **@{-}, (10,88);(10,87) **@{-},
(15,88);(15,87) **@{-}, (20,88);(20,87) **@{-},
(25,88);(25,87) **@{-}, (30,88);(30,87) **@{-},
(35,88);(35,87) **@{-}, (40,88);(40,87) **@{-},
(45,88);(45,87) **@{-}, (50,88);(50,87) **@{-},
(55,88);(55,87) **@{-}, (60,88);(60,87) **@{-},
(65,88);(65,87) **@{-}, (70,88);(70,87) **@{-},
(75,88);(75,87) **@{-}, (80,88);(80,87) **@{-},
(85,88);(85,87) **@{-},
(5,83)*{_{\text{i}}}, (10,83)*{_{+}},
(5,78)*{_{\text{i}}}, (10,78)*{_{\text{i}}}, (15,78)*{_{+}},
(5,74)*{_{\vdots}}, (10,73)*{_{\text{i}}}, (15,73)*{_{\text{i}}}, (20,74)*{_{\ddots}},
(5,69)*{_{\vdots}}, (10,69)*{_{\vdots}}, (15,68)*{_{\text{i}}}, (20,69)*{_{\ddots}}, (25,69)*{_{\ddots}},
(5,64)*{_{\vdots}}, (10,64)*{_{\vdots}}, (15,64)*{_{\vdots}}, (20,64)*{_{\ddots}}, (25,64)*{_{\ddots}}, (30,64)*{_{\ddots}},
(5,58)*{_{\text{i}}}, (10,58)*{_{\text{i}}}, (15,58)*{_{\cdots}}, (20,58)*{_{\text{i}}}, (25,58)*{_{\text{i}}}, (30,58)*{_{\text{i}}}, (35,58)*{_{+}},
(5,53)*{_{\text{i}}}, (10,53)*{_{\text{i}}}, (15,53)*{_{\cdots}}, (20,53)*{_{\cdots}}, (25,53)*{_{\text{i}}}, (30,53)*{_{\text{i}}}, (35,53)*{_{\text{i}}}, (40,53)*{_{+}},
(5,48)*{{\diamond}}, (10,48)*{_{\text{I}}}, (15,48)*{_{\text{I}}}, (20,48)*{_{\cdots}}, (25,48)*{_{\cdots}}, (30,48)*{_{\text{I}}}, (35,48)*{_{\text{I}}}, (40,48)*{_{0}}, (45,48)*{_{+}},
(5,43)*{{\diamond}}, (10,43)*{{\diamond}}, (15,43)*{_{\text{I}}}, (20,43)*{_{\text{I}}}, (25,43)*{_{\cdots}}, (30,43)*{_{\cdots}}, (35,43)*{_{\text{I}}}, (40,43)*{_{1}}, (45,43)*{_{0}}, (50,43)*{_{+}},
(5,38)*{\diamond}, (10,38)*{\diamond}, (15,38)*{\diamond}, (20,38)*{_{\text{I}}}, (25,39)*{_{\ddots}}, (30,39)*{_{\ddots}}, (35,39)*{_{\vdots}}, (40,38)*{_{2}}, (45,38)*{_{1}}, (50,38)*{_{0}}, (55,38)*{_{+}},
(5,34)*{_{\vdots}}, (10,34)*{_{\ddots}}, (15,34)*{_{\ddots}}, (20,34)*{_{\ddots}}, (25,34)*{_{\ddots}}, (30,34)*{_{\ddots}}, (35,34)*{_{\vdots}}, (40,33)*{_{3}}, (45,33)*{_{2}}, (50,33)*{_{1}}, (55,33)*{_{0}}, (60,33)*{_{+}},
(5,29)*{_{\vdots}},(10,29)*{_{\vdots}},(15,29)*{_{\vdots}},(20,29)*{_{\vdots}},(25,29)*{_{\ddots}}, (30,29)*{_{\ddots}}, (35,29)*{_{\vdots}}, (40,29)*{_{\vdots}}, (45,29)*{_{\ddots}}, (50,29)*{_{\ddots}}, (55,29)*{_{\ddots}}, (60,29)*{_{\ddots}}, (65,29)*{_{\ddots}},
(25,24)*{_{\vdots}},(30,23)*{\diamond}, (35,23)*{_{\text{I}}}, (40,24)*{_{\vdots}}, (45,24)*{_{\ddots}}, (50,23)*{_{3}}, (55,23)*{_{2}}, (60,23)*{_{1}}, (65,23)*{_{0}}, (70,23)*{_{+}},
(30,19)*{_{\vdots}},(35,18)*{\diamond}, (40,18)*{_{_{m_0-1}}}, (45,18)*{_{\cdots}}, (50,18)*{_{\cdots}}, (55,18)*{_{3}}, (60,18)*{_{2}}, (65,18)*{_{1}}, (70,18)*{_{0}}, (75,18)*{_{+}},
(35,14)*{_{\vdots}},(40,13)*{\diamond}, (45,13)*{_{_{m_0-1}}}, (50,13)*{_{\cdots}}, (55,13)*{_{\cdots}}, (60,13)*{_{3}}, (65,13)*{_{2}}, (70,13)*{_{1}}, (75,13)*{_{0}}, (80,13)*{_{+}},
(40,9)*{_{\vdots}},(45,8)*{\diamond}, (50,8)*{_{_{m_0-1}}}, (55,8)*{_{\cdots}}, (60,8)*{_{\cdots}}, (65,8)*{_{3}}, (70,8)*{_{2}}, (75,8)*{_{1}}, (80,8)*{_{0}}, (85,9)*{_{\ddots}},
(45,4)*{_{\vdots}},(50,4)*{_{\ddots}}, (55,4)*{_{\ddots}}, (60,4)*{_{\ddots}}, (65,4)*{_{\ddots}}, (70,4)*{_{\ddots}}, (75,4)*{_{\ddots}}, (80,4)*{_{\ddots}}, (85,4)*{_{\ddots}},
(6,14);(14,6) **@{-}, (14,6);(22,14)**@{-}, (22,14);(14,22) **@{-}, (14,22);(6,14)**@{-},
(2.5,88);(88,2.5) **@{.}, (2.5,88);(2.5,2.5) **@{.}, (2.5,50.5);(40,50.5) **@{.}, (37,50.5);(37,18) **@{.}, (4.5,50.5);(55,0) **@{.},
(-2,90)*{_{0}},
(-3,83)*{_{1}}, (-3,78)*{_{2}},
(-3,73)*{_{3}}, (-3,69)*{_{\vdots}},
(-3,64)*{_{\vdots}}, (-5,58)*{_{m_0-1}},
(-3,53)*{_{m_0}}, (-5,48)*{_{m_0+1}},
(-5,43)*{_{m_0+2}}, (-5,38)*{_{m_0+3}},
(-3,34)*{_{\vdots}}, (-3,29)*{_{\vdots}},
(-5,23)*{_{2m_0-1}}, (-3,18)*{_{2m_0}},
(-5,13)*{_{2m_0+1}}, (-3,9)*{_{\vdots}},
(-3,4)*{_{\vdots}}, (-8,3)*{{n}},
(5,91)*{_{1}}, (10,91)*{_{2}},
(15,91)*{_{3}}, (20,91)*{_{\cdots}},
(25,91)*{_{\cdots}}, (30,91)*{_{\cdots}},
(35,91)*{_{_{m_0}}}, (40,91)*{_{}},
(45,91)*{_{_{m_0+2}}}, (50,91)*{_{}},
(52,91)*{_{\cdots}}, (58,91)*{_{\cdots}},
(64,91)*{_{\cdots}}, (70,91)*{_{_{2m_0}}},
(75,91)*{_{}}, (80,91)*{_{_{2m_0+2}}},
(87,91)*{_{\cdots}}, (91,92)*{i},
(60,80)*{{0,...,m_0-1 : \text{diagonals per Prop. 3.2 (i)}}},
(70,75)*{{\text{i : check, per Prop. 3.2 (ii)}}},
(70,70)*{{\text{I : check, per Prop. 3.2 (iii)}}},
(70,65)*{{+ : \text{automatically positive}}},
(70,60)*{{\diamond : \text{sign as term to above}}},
(45,-5)*{\text{Figure 1}}
\endxy} \vspace*{5pt}
\end{figure}

\bigskip

In Figure 1, the symbol ``+'' denotes a term automatically
positive (such terms on the coordinate axes are omitted); ``i'' (respectively, ``I'')
denotes terms that must be checked individually as in condition (ii)
(respectively, condition (iii)) of the theorem; ``$\diamond$%
'' denotes terms guaranteed non-negative if the term immediately
above is non-negative; the numbered diagonals correspond to the $c(n,n-i)$
as in condition (i) of the theorem, for the various values of $i$, $0 \leq i
\leq m_0$.

\bigskip

\textbf{4. Main theorem.} For the reader's convenience, we recall Stampfli's
subnormal completion (cf. \cite{CuF2}, \cite{St}). For given numbers $\alpha
_{0},\alpha _{1},\alpha _{2}$ with $0<\alpha _{0}<\alpha _{1}<\alpha _{2}$,
define
\begin{equation}
\alpha _{n}^{2}=\Psi _{1}+\frac{\Psi _{0}}{\alpha _{n-1}^{2}}\ \ \text{for
all}\ n\geq 3,  \tag{4.1}
\end{equation}%
where $\Psi _{0}=-\frac{\alpha _{0}^{2}\alpha _{1}^{2}\left( \alpha
_{2}^{2}-\alpha _{1}^{2}\right) }{\alpha _{1}^{2}-\alpha _{0}^{2}}$ and $%
\Psi _{1}=\frac{\alpha _{1}^{2}\left( \alpha _{2}^{2}-\alpha _{0}^{2}\right)
}{\alpha _{1}^{2}-\alpha _{0}^{2}}$. Then we may obtain a weight sequence $%
\{\alpha _{n}\}_{n=0}^{\infty }$ generated recursively by (4.1), which is
usually denoted by $(\alpha _{0},\alpha _{1},\alpha _{2})^{\wedge }$, and the resulting shift is subnormal (see \cite{St}). It follows from \cite{CuF2} that
\begin{equation*}
\alpha _{n}\nearrow L:=\frac{1}{\sqrt{2}}\left( \Psi _{1}+\sqrt{\Psi
_{1}^{2}+4\Psi _{0}}\right) ^{1/2}\ \ \text{as}\ \ n\rightarrow \infty .
\end{equation*}

Let $\alpha =\{\alpha _{n}\}_{n=0}^{\infty }=\sqrt{x_{m}},\sqrt{x_{m-1}}%
,\cdots ,\sqrt{x_{1}},(\sqrt{u},\sqrt{v},\sqrt{w})^{\wedge }$ with $%
0<x_{m}\leq \cdots \leq x_{1}\leq u<v<w$ be a sequence which is called
a \emph{backward} $m$-\textit{step extension} of $(\sqrt{u},\sqrt{v},\sqrt{w}%
)^{\wedge }$. For brevity, we denote such a backward $m$-step extension by $%
\widetilde{\alpha }^{(m)}(u,v,w)$.  Recall that we define $\eta_n = \frac{v_n}{u_n}$ where $u_n$ and $v_n$ are as defined before (2.2).

\medskip

\textbf{Remark.}  Before continuing we must dispose of a certain technical obstruction.  Consider the backward extension with weight sequence
$$\sqrt{1}, \sqrt{2},\sqrt{3},\sqrt{3},\sqrt{3},\sqrt{3},\sqrt{3},\sqrt{3},\sqrt{5},
\sqrt{6},(\sqrt{u},\sqrt{v},\sqrt{w})^{\wedge },$$
(where of course $6 \leq u < v < w$ by our standing assumption).  The presence of the six equal weights produces the unwelcome fact that $v_5 = 0$, which violates the standing assumption of Section 3, means that all but finitely many of the $c(n,n+1)$ are zero, and casts into doubt the strict positivity of determinants of some submatrices.  Remark 2.4 from \cite{CuF1} shows that this strict positivity may not be done without.  We leave to the reader to verify (but see Proposition 2.8 of \cite{CuF2}) that in this case the matrix $D_n$ splits into a direct sum of a finite matrix of the type we have been considering, a zero matrix, and (the finite restriction of) an infinite matrix corresponding to a shorter backward extension of the kind we are considering.  More blocks of six or more successive equal weights in the backward extension, of course, create direct sums with more terms.  In these cases the positivity of the finite matrices must be analyzed separately, as well as using our analysis to analyze the (restrictions of) the infinite portion.  For ease of exposition we will omit discussion of the analysis of these finite matrices, and assume without further comment that we are without this difficulty in what follows. Note finally that the examples in Section 5 certainly meet this assumption.

\medskip

\textbf{Lemma 4.1.} \textit{Let }$\widetilde{\alpha }^{(m)}(u,v,w)$\textit{\
be a backward }$m$\textit{-step extension. Then }$u_{n}v_{n+2}=w_{n}\
(n\geq m+1).$

\smallskip

\textit{Proof}. Imitate the proof of \cite[Lemma 3.1]{LLB}. \hfill$\Box$

\medskip

\textbf{Lemma 4.2} \textit{Let }$\widetilde{\alpha }^{(m)}(u,v,w)$\textit{\
be a backward }$m$\textit{-step extension. Then}

(i) $\eta _{n+1}\geq \eta _{n}$ \textit{for all }$n\geq m+3$,

(ii) $\lim_{n\rightarrow \infty }\eta _{n}=K:=\frac{(\Psi _{1}^{2}+\Psi
_{0})^{2}}{\Psi _{0}^{2}}L^{4}$,

\noindent \textit{where} $L=\frac{1}{\sqrt{2}}\left( \Psi _{1}+\sqrt{\Psi _{1}^{2}+4\Psi _{0}}%
\right) ^{1/2}$.

\smallskip

\textit{Proof}. Imitate the proof of \cite[Lemma 3.3]{LLB}. \hfill$\Box$

\medskip

We denote
\begin{align*}
d^{[1]}(t)&:=\text{det}D_{n}^{[1]}(t)=\sum_{i=0}^{n+1}c^{[1]}(n,i)t^{i}, \\
d^{[2]}(t)&:=\text{det}D_{n}^{[2]}(t)=\sum_{i=0}^{n+1}c^{[2]}(n,i)t^{i}.
\end{align*}

The task at hand is to apply the machinery of Section 3 to the two (families of) matrices $D_{n}^{[1]}(t)$ and $D_{n}^{[2]}(t)$.  This will involve replacing various objects $\check{\square}$ in that section with their $\square^{[1]}$ and $\square^{[2]}$ versions.  We will use expressions like $G^{[j]}$, $Q^{[j]}$, and $\eta^{[j]}$ with the obvious meanings.  It may be useful to note in advance that the limit $\lim_{n\rightarrow \infty }\check{\eta}_{n}=\check{\eta} <\infty$ turns out to be the same for both ``$[1]$'' and ``$[2]$'' (see Lemma 4.2 above), as do certain quantities $A$ and $B$.  To ease the burden of notation slightly we omit the $[j]$ superscripts on these quantities.

To begin, taking $m_{0}=\left[ \frac{m}{2}\right] +1$, where $[x]$ is the Gauss number
of $x$, we obtain
\begin{align}
u_{2n}v_{2(n+1)}& =w_{2n},\ \ u_{2n+1}v_{2(n+1)+1}=w_{2n+1},\ n\geq m_{0},
\tag{4.1} \\
\eta _{2(n+1)}& \geq \eta _{2n},\ \ \eta _{2(n+1)+1}\geq \eta _{2n+1},\
n\geq m_{0}+1.  \notag
\end{align}%
Then, if we apply Proposition 3.1 and Proposition 3.2 using  (4.1) to $%
D_{n}^{[1]}$ and $D_{n}^{[2]}$ as in (2.2), Proposition 3.1 and Proposition 3.2
can be restated as follows.

\textbf{Proposition 4.3.} \textit{Let }$\widetilde{\alpha }^{(m)}(u,v,w)$%
\textit{\ be a backward }$m$\textit{-step extension and }$m_{0}=%
\left[ \frac{m}{2}\right] +1$. \textit{Suppose that for each }$i$, $\check{G}(n,i)$ \textit{is monotone increasing or decreasing for }$n\geq m_0 + i +1$ \textit{for both }$j=1$ \textit{and } $j=2$.  \textit{For either }$j=1$ \textit{or } $j=2$, \textit{the following conditions are equivalent:}

(i) $c^{[j]}(n,n-i)\geq 0$ (\textit{all} $n\geq m_{0}+i+1 $),

(ii)\textit{\ it holds that either}

\ \ \ (ii-a) $G^{[j]}(n,i)$ \textit{is increasing in} $n$, $n\geq m_{0}+i+1$
\textit{and} $G^{[j]}(m_{0}+1+i,i)\geq 0$,

\textit{or}

\ \ \ (ii-b) $G^{[j]}(n,i)$ \textit{is decreasing in} $n$, $n\geq m_{0}+i+1$
\textit{and}
\begin{equation*}
v^{[j]}_{m_{0}}\cdots v^{[j]}_{0},+\sum_{k=0}^{i}K ^{k+1}\widehat{c}%
^{[j]}(m_{0},m_{0}-k)\geq 0.
\end{equation*}
\emph{where } $\widehat{c}^{[j]}(n,i)=v^{[j]}_n c^{[j]}(n-1,i-1)-w^{[j]}%
_{n-1}c^{[j]}(n-2,i-1)$, $j=1,2$.

\medskip

\textbf{Proposition 4.4.} \textit{Let }$\widetilde{\alpha }^{(m)}(u,v,w)$%
\textit{\ be a backward }$m$\textit{-step extension and }$m_{0}=\left[ \frac{%
m}{2}\right] +1$. \textit{Then }$D_{n}^{[j]},$%
\textit{\ $j=1$ or $j=2$,} \textit{\ has
positive determinant coefficients if and only if }

(i) \textit{for each }$i$\textit{, }$0\leq i\leq m_{0}-1$\textit{, one of
the conditions of Proposition 4.3 holds,}

(ii) \textit{for each }$n$\textit{, }$1\leq n\leq m_{0}$\textit{\ and each }$%
i$\textit{, }$1\leq i\leq n$\textit{, }$c^{[j]}(n,i)\geq 0$,

(iii) \textit{for each }$n$\textit{, }$m_{0}+1\leq n\leq 2m_{0}-1$\textit{\
and each }$i$\textit{, }$n+1-m_{0}\leq i\leq m_{0}$\textit{, }$%
c^{[j]}(n,i)\geq 0$\textit{.}

\medskip

The statements of Proposition 4.4 can be improved by using the following
proposition.

\medskip

\textbf{Proposition 4.5} \textit{Let }$\widetilde{\alpha }^{(m)}(u,v,w)$%
\textit{\ be a backward }$m$\textit{-step extension and }$m_{0}=\left[ \frac{%
m}{2}\right] +1$. \emph{Then for either $j=1$ or $j=2$, $Q^{[j]}(\ell,k)$ is constant in} $\ell$, $\ell \geq m_0 +1$.

\medskip

Before giving the proof of Proposition 4.5, we first consider some crucial
lemmas below.

Note that if $(\sqrt{u},\sqrt{v},\sqrt{w})^{\wedge }=\sqrt{u},(\sqrt{v},%
\sqrt{w},\sqrt{p})^{\wedge }$ for given $u,v,w,$ by simple computation, we
have
\begin{equation}
p=\frac{v(w^{2}+vu-2wu)}{w(v-u)}.  \tag{4.2}
\end{equation}
Then a direct computation provides the following lemma.

\medskip

\textbf{Lemma 4.6.} \textit{Let}
\begin{align*}
A(u,v,w)& :=\frac{v^{3}(uv(v-3w)+u^{2}w+vw^{2})^{6}}{%
u^{2}(u-v)^{8}(v-w)^{2}(u^{2}(v-2w)-2uv(v-2w)-vw^{2})}, \\
B(u,v,w)& :=\frac{%
v(uv(v-3w)+u^{2}w+vw^{2})^{3}(u^{2}(2v-3w)-vw^{2}+uv(-3v+5w))}{%
u^{2}(u-v)^{4}(v-w)^{2}(u^{2}(v-2w)-2uv(v-2w)-vw^{2})},
\end{align*}%
\textit{and suppose that }$(\sqrt{u},\sqrt{v},\sqrt{w})^{\wedge }=\sqrt{u},(%
\sqrt{v},\sqrt{w},\sqrt{p})^{\wedge }$ \textit{for} $0<u<v<w<p.$ \textit{Then%
}
\begin{equation*}
A(u,v,w)=A(v,w,p)\mathit{\ }\text{\textit{and }}B(u,v,w)=B(v,w,p).
\end{equation*}

\medskip

For brevity, we let $A:=A(u,v,w)$ and $B:=B(u,v,w)$.

\bigskip

\textbf{Lemma 4.7.} \textit{A Stampfli completion }$(\sqrt{u},\sqrt{v},%
\sqrt{w})^{\wedge }$\textit{\ satisfies}
\begin{equation*}
\eta _{n}=\frac{A}{\eta _{n-2}\eta _{n-4}}+B,n\geq 7.
\end{equation*}%

\smallskip

\textit{Proof.} For mathematical induction, we consider
\begin{equation*}
\eta _{7}=\frac{A}{\eta _{5}\eta _{3}}+B,
\end{equation*}%
which follows from an easy computation. And also observe that (in notation whose meaning is obvious)
\begin{align*}
\eta _{8}((\sqrt{u},\sqrt{v},\sqrt{w})^{\wedge })& =\eta _{8}(\sqrt{u},(%
\sqrt{v},\sqrt{w},\sqrt{p})^{\wedge }) \\
& =\eta _{7}((\sqrt{v},\sqrt{w},\sqrt{p})^{\wedge }) \\
& =\frac{A(v,w,p)}{\eta _{5}((\sqrt{v},\sqrt{w},\sqrt{p})^{\wedge })\eta
_{3}((\sqrt{v},\sqrt{w},\sqrt{p})^{\wedge })}+B(v,w,p) \\
& =\frac{A}{\eta _{6}((\sqrt{u},\sqrt{v},\sqrt{w})^{\wedge })\eta _{4}((%
\sqrt{u},\sqrt{v},\sqrt{w})^{\wedge })}+B
\end{align*}%
and%
\begin{align*}
\eta _{9}((\sqrt{u},\sqrt{v},\sqrt{w})^{\wedge })& =\eta _{8}((\sqrt{v},%
\sqrt{w},\sqrt{p})^{\wedge }) \\
& =\frac{A(v,w,p)}{\eta _{6}((\sqrt{v},\sqrt{w},\sqrt{p})^{\wedge })\eta
_{4}((\sqrt{v},\sqrt{w},\sqrt{p})^{\wedge })}+B(v,w,p) \\
& =\frac{A}{\eta _{7}((\sqrt{u},\sqrt{v},\sqrt{w})^{\wedge })\eta _{5}((%
\sqrt{u},\sqrt{v},\sqrt{w})^{\wedge })}+B.
\end{align*}%
Repeating the argument proves this lemma. \hfill$\Box$

\bigskip

For a backward $m$-step extension $\widetilde{\alpha }^{(m)}(u,v,w),$ Lemma
4.7 can be restated as follows.

\medskip

\textbf{Lemma 4.8.} \textit{Let }$\widetilde{\alpha }^{(m)}(u,v,w)$ \textit{%
be as usual and }$m_{0}=\left[ \frac{m}{2}\right] +1$. \textit{Let }$\check{\eta%
}_{n}=\eta _{2n}$\textit{\ }$[$\textit{or }$\check{\eta}_{n}=\eta _{2n+1}].$%
\textit{\ Then}
\begin{equation*}
\check{\eta}_{n+2}=\frac{A}{\check{\eta}_{n+1}\check{\eta}_{n}}+B,n\geq
m_{0}+1.
\end{equation*}

\medskip

In the next two proofs we adopt temporarily that convention that ``$\check{\square}$'' will stand for ``$\square^{[j]}$'' for any quantity ``$\square$'' as convenient (of course $j=1$ or $j=2$).  Observe that by this convention the previous lemma is a statement about either the $\eta^{[1]}_n = \eta_{2n}$ or $\eta^{[2]}_n = \eta_{2n+1}$ sequence as needed. With this convention in mind we will also use $j$ temporarily as a running subscript.

\medskip

\textbf{Lemma 4.9.} \textit{Let }$\widetilde{\alpha }^{(m)}(u,v,w)$ \textit{%
be as usual and }$m_{0}=\left[ \frac{m}{2}\right] +1$.\textit{\ Then for }$%
\ell \geq m_0 +1$,
\begin{equation*}
\check{Q}(\ell,k+2)=\left\{
\begin{array}{ll}
A\check{Q}(\ell,k-1)+B\check{Q}(\ell,k+1), & k\geq 1, \\
B\check{Q}(\ell,1), & k=0.%
\end{array}%
\right.
\end{equation*}

\smallskip

\textit{Proof}. The result is easy for $k=0$.  Using Lemma 4.8, we have
\begin{align*}
\ \ \ \ \ \ \ \ \ \ \check{Q}(\ell,k+2)& =\frac{\prod_{j=0}^{k+2}\check{\eta}_{\ell+j+1}-\prod_{j=-1}^{k+1}%
\check{\eta}_{\ell+j+1}}{\check{\eta}_{\ell+1}-\check{\eta}_{\ell}} \\
& =\frac{\check{\eta}_{\ell+k+3}\prod_{j=0}^{k+1}\check{\eta}_{\ell+j+1}-%
\check{\eta}_{\ell+k+2}\prod_{j=-1}^{k}\check{\eta}_{\ell+j+1}}{\check{\eta}%
_{\ell+1}-\check{\eta}_{\ell}} \\
& =\frac{\left( \frac{A}{\check{\eta}_{\ell+k+2}\check{\eta}_{\ell+k+1}}%
+B\right) \prod_{j=0}^{k+1}\check{\eta}_{\ell+j+1}}{\check{\eta}_{\ell+1}-%
\check{\eta}_{n-i}} \\
& \text{ \ \ \ \ \ \ \ \ \ \ \ \ \ \ \ }-\frac{\left( \frac{A}{\check{\eta}%
_{\ell+k+1}\check{\eta}_{\ell+k}}+B\right) \prod_{j=-1}^{k}\check{\eta}%
_{\ell+j+1}}{\check{\eta}_{\ell+1}-\check{\eta}_{\ell}} \\
& =A\left( \frac{\prod_{j=0}^{k-1}\check{\eta}_{\ell+j+1}-\prod_{j=-1}^{k-2}%
\check{\eta}_{\ell+j+1}}{\check{\eta}_{\ell+1}-\check{\eta}_{\ell}}\right) \\
& \text{ \ \ \ \ \ \ \ \ \ \ }+B\left( \frac{\prod_{j=0}^{k+1}\check{\eta}%
_{\ell+j+1}-\prod_{j=-1}^{k}\check{\eta}_{\ell+j+1}}{\check{\eta}_{\ell+1}-%
\check{\eta}_{\ell}}\right) \\
& =A\check{Q}(\ell,k-1)+B\check{Q}(\ell,k+1).\ \ \ \ \ \ \ \ \ \ \ \ \ \ \ \ \ \ \ \ \ \ \ \ \ \ \ \ \ \ \hfill\Box
\end{align*}

\medskip

We now give the proof of Proposition 4.5 below.

\medskip

\textit{Proof of Proposition 4.5.} It follows obviously from (3.3) that
\begin{equation*}
\check{Q}(\ell,0)=\frac{\check{\eta}_{\ell+1}-\check{\eta}_{\ell}}{\check{\eta}_{\ell+1}-%
\check{\eta}_{\ell}}=1.
\end{equation*}%
Using the idea of Lemma 4.4 in \cite{EJP} (cf. \cite[Lemma 3.2]{BEJL}), we
have
\begin{equation*}
\frac{\eta _{n}\eta _{n+2}-\eta _{n-2}\eta _{n}}{\eta _{n}-\eta _{n-2}}
\end{equation*}%
is constant in $n$ for $n\geq m+5$. Then it follows that
\begin{equation*}
\check{Q}(\ell,1)=\frac{\check{\eta}_{\ell+1}\check{\eta}_{\ell+2}-\check{\eta}_{\ell}%
\check{\eta}_{\ell+1}}{\check{\eta}_{\ell+1}-\check{\eta}_{\ell}}
\end{equation*}%
is constant in $\ell$, for $\ell \geq m_{0}+1$. Hence by Lemma 4.9,
\begin{equation*}
\check{Q}(\ell,k)=\check{\Gamma} _{0}(A,B)+\check{\Gamma }_{1}(A,B)\check{Q}(\ell,1)
\end{equation*}%
is constant in $\ell$, for $\ell\geq m_{0}+1$, for some constants $\check{\Gamma}
_{0}(A,B)$ and $\check{\Gamma} _{1}(A,B)$ independent of $\ell$. Hence the proof is complete. $%
\hfill \Box $

\medskip

Since (interpreting the result in the ``$[j]$'' notation), $Q^{[j]}(\ell,k)$ is constant in $\ell$ for $\ell \geq m_0 + 1$, we may write $Q^{[j]}_{k}\equiv
Q^{[j]}_{k}(u,v,w):=Q(\ell,k)$, for $\ell \geq m_{0}+1$. It will be useful to note in what follows that this is applied with $\ell$ set to $n-i$, including implicitly in part (ii-b) of the next proposition.  According to Proposition 4.5,
Proposition 4.3 can be restated usefully as follows.

\medskip

\textbf{Proposition 4.10.} \textit{With notation as above, fix} $i$ \textit{%
such that} $0\leq i\leq m_{0}-1$, \textit{and let} $j=1$ \textit{or } $j=2$. \textit{Then the
following conditions are equivalent:}

(i) \textit{the coefficients satisfy} $c^{[j]}(n,n-i)\geq 0$ (\textit{all} $n\geq
m_{0}+i+1 $),

(ii)\textit{\ it holds that either}

\ \ \ \ (ii-a) $\sum_{k=0}^{i}Q^{[j]}_{k}\widehat{c}^{[j]}(m_{0},m_{0}-k)\geq 0$
\textit{and} $c^{[j]}(m_{0}+1+i,m_{0}+1)\geq 0$, \textit{or}

\ \ \ \ (ii-b) $\sum_{k=0}^{i}Q^{[j]}_{k}\widehat{c}^{[j]}(m_{0},m_{0}-k)\leq 0$
\textit{and}
\begin{equation*}
\lim_{n\rightarrow \infty }G^{[j]}(n,i)=\check{v}_{m_{0}}\cdots \check{v}%
_{0}+\sum_{k=0}^{i}K^{k+1}\widehat{c}^{[j]}(m_{0},m_{0}-k)\geq 0.
\end{equation*}

\medskip

Hence the following theorem comes from Propositions 4.4 and 4.10.

\medskip

\textbf{Theorem 4.11. }\textit{Let }$\widetilde{\alpha }^{(m)}:=\widetilde{%
\alpha }^{(m)}(u,v,w)$\textit{\ be a backward} $m$\textit{-step extension
and }$m_{0}=\left[ \frac{m}{2}\right] +1$.\textit{\ Then }$W_{\widetilde{%
\alpha }^{(m)}}$\textit{\ is} \textit{a semi-cubically hyponormal weighted
shift with p.d.c. if and only if the following conditions hold for each} $%
j=1,2$\textit{:}

(i) \textit{for each }$i$\textit{, }$0\leq i\leq m_{0}-1$\textit{, one of
the second conditions of Proposition 4.10 holds,}

(ii) \textit{for each }$n$\textit{, }$1\leq n\leq m_{0}$\textit{\ and each }$%
i$\textit{, }$1\leq i\leq n$\textit{, }$c^{[j]}(n,i)\geq 0$\textit{,}

(iii) \textit{for each }$n$\textit{, }$m_{0}+1\leq n\leq 2m_{0}-1$\textit{\
and each }$i$\textit{, }$n+1-m_{0}\leq i\leq m_{0}$\textit{, }$%
c^{[j]}(n,i)\geq 0$\textit{.}

\bigskip

\textbf{5. Examples.} In 1991, R. Curto (\cite{Cu2}) suggested a problem:
``describe all quadratically hyponormal weighted shifts $%
W_{\alpha }$ with first two weights equal.'' Subsequently
several papers have considered the problem of quadratic
hyponormality of $W_{\alpha }$ with first two weights equal (cf. \cite{CuF1}%
, \cite{CuF2}, \cite{CuJ}, \cite{DEJL}, \cite{EJP}, \cite{LLB}). On the
other hand, in \cite{LCL}, it was shown that if $W_{\alpha }$ is cubically
hyponormal with first two weights equal, then $W_{\alpha }$ has flatness.
And, in \cite{DEJL}, it was proved that there exists a semi-cubically
hyponormal weighted shift $W_{\alpha }$ with $\alpha _{0}=\alpha _{1}<\alpha
_{2}$ but which is not 2-hyponormal. These facts motivate a parallel problem
``describe all semi-cubically hyponormal weighted shifts $%
W_{\alpha }$ with first two weights equal.'' As part of this study, Do-Exner-Jung-Li(\cite{DEJL}) characterized the semi-cubic
hyponormality of the weighted shift $W_{\alpha (x)}$ with positive
determinant coefficients (p.d.c.), where $\alpha (x):\sqrt{x},\sqrt{x},\sqrt{%
\frac{k+1}{k+2}}$ $(k\geq 2)$ is a weight sequence with Bergmann tail. In
\cite{BEJL} one also described the semi-cubic hyponormality
of the weighted shift $W_{\alpha }$ with the condition p.d.c., which is
associated to a weight sequence $\alpha :1,1,\sqrt{x},(\sqrt{u},\sqrt{v},%
\sqrt{w})^{\wedge }$. Note that the conditions of the case $m=3$ in Theorem
4.11 coincide with the conditions in \cite[Theorem 3.4]{BEJL}. As a
continued study, we consider a backward 4-step extension $\alpha :1,1,%
\sqrt{x},\sqrt{y},\left( \sqrt{u},\sqrt{v},\sqrt{w}\right) ^{\wedge }$ with $%
1<x<y<u<v<w$. According to Theorem 4.11, we describe the semi-cubic
hyponormality of such weighted shifts $W_{\alpha }$ below.  To ease the burden of notation, in what follows we omit (leave implicit) ``$[j]$'' on instances of $Q$ which appear in an expression already containing this superscript;  in no expression will ``$[1]$'' and ``$[2]$'' be mingled.

\medskip

According to Theorem 4.11, for each $j=1,2$, $D_{n}^{[j]}(t)$ has positive
determinants coefficients for all $n$ if and only if the following
conditions hold:

\smallskip

\noindent (i) $c^{[j]}(1,1)$, $c^{[j]}(2,1)$, $c^{[j]}(3,1)$, $c^{[j]}(2,2)$%
, $c^{[j]}(3,2)$, $c^{[j]}(4,2)$, $c^{[j]}(3,3)$, $c^{[j]}(4,3)$, $%
c^{[j]}(5,3) $ are positive,

\smallskip

\noindent (ii) one of following conditions holds

\smallskip

(ii-a) $\widehat{c}^{[j]}(3,3)\geq 0$, $c^{[j]}(4,4)\geq 0$ or,

\smallskip

(ii-b) $\widehat{c}^{[j]}(3,3)\leq 0$, $v^{[j]}_3 v^{[j]}_2 v^{[j]}_1 v^{[j]}_0+K\widehat{c}^{[j]}(3,3)\geq 0$,

\smallskip

\noindent (iii) one of following conditions holds;

\smallskip

(iii-a) $\sum_{k=0}^{1}Q_{k}\widehat{c}^{[j]}(3,3-k)\geq 0$, $%
c^{[j]}(5,4)\geq 0$ or,

\smallskip

(iii-b) $\sum_{k=0}^{1}Q_{k}\widehat{c}^{[j]}(3,3-k)\leq 0$, $v^{[j]}_3 v^{[j]}_2 v^{[j]}_1 v^{[j]}_0+\sum_{k=0}^{1}K^{k+1}\widehat{c}%
^{[j]}(3,3-k)\geq 0$,

\smallskip

\noindent (iv) one of following conditions holds

\smallskip

(iv-a) $\sum_{k=0}^{2}Q_{k}\widehat{c}^{[j]}(3,3-k)\geq 0$, $%
c^{[j]}(6,4)\geq 0$ or,

\smallskip

(iv-b) $\sum_{k=0}^{2}Q_{k}\widehat{c}^{[j]}(3,3-k)\leq 0$, $v^{[j]}_3 v^{[j]}_2 v^{[j]}_1 v^{[j]}_0+\sum_{k=0}^{2}K^{k+1}\widehat{c}%
^{[j]}(3,3-k)\geq 0,$

\smallskip

\noindent where $\widehat{c}^{[j]}(3,3-k)=v^{[j]}_3 c^{[j]}(2,2-k)-w^{[j]}_{2}c^{[j]}(1,2-k)$, and $K$ is as in Lemma 4.2.

\medskip

In preparation for what follows, we record some computations, beginning with
\begin{eqnarray*}
c^{[1]}(1,1) &=&(uy-1)x>0, \\
c^{[2]}(2,1) &=&(v-u)c^{[2]}(1,1), \\
c^{[2]}(3,1) &=&\frac{u(w-v)^{2}}{w}c^{[2]}(1,1).
\end{eqnarray*}%
One computes that
\begin{equation*}
c^{[2]}(1,1)>0\Longleftrightarrow c^{[2]}(2,1)>0\Longleftrightarrow
c^{[2]}(3,1)>0.
\end{equation*}%
By some further computations (using, when convenient, the trick of setting $v$ to $u + \delta$ and $w$ to $u + \delta + \epsilon$ with both $\delta$ and $\epsilon$ positive)
\begin{equation*}
c^{[2]}(2,2)>0\Longleftrightarrow c^{[2]}(3,2)>0\Longleftrightarrow
c^{[2]}(4,2)>0,
\end{equation*}%
\begin{equation*}
c^{[2]}(3,3)>0\Longleftrightarrow c^{[2]}(4,3)>0\Longleftrightarrow
c^{[2]}(5,3)>0,
\end{equation*}%
and
\begin{equation*}
c^{[2]}(4,4)>0\Longleftrightarrow c^{[2]}(5,4)>0\Longleftrightarrow
c^{[2]}(6,4)>0.
\end{equation*}%
Observe that $\widehat{c}^{[2]}(3,2)=\widehat{c}^{[2]}(3,1)=0$, and by Proposition 3.9, we have (for $j = 1, 2$)
\begin{align*}
Q^{[j]}_{0}& =1, \\
Q^{[j]}_{1}& =\frac{%
v(2uv^{2}+2u^{2}w+vw^{2}-u^{2}v-4uvw)(uv^{2}+u^{2}w-3uvw+vw^{2})^{2}}{%
u^{2}(v-u)^{4}(w-v)^{2}}, \\
Q^{[j]}_{2}& =BQ^{[j]}_{1}.
\end{align*}

\noindent Using the conditions above, we may characterize the semi-cubic
hyponormality of a weighted shift $W_{\alpha }$ with weight sequence $\alpha :1,1,%
\sqrt{x}, \sqrt{y}, \left( \sqrt{u}, \sqrt{v}, \sqrt{w} \right) ^{\wedge }$ with $%
1<x<y<u<v<w$ and with the p.d.c. condition as follows.

\medskip

\textbf{Proposition 5.1.}\textit{\ Let }$\alpha :1,1,\sqrt{x},\sqrt{y}%
,\left( \sqrt{u},\sqrt{v},\sqrt{w}\right) ^{\wedge }$\textit{\ with }$%
1<x<y<u<v<w$\textit{. Then }$W_{\alpha }$\textit{\ is semi-cubically
hyponormal weighted shift with p.d.c. if and only if}

\smallskip

(i) $c^{[1]}(2,1)$, $c^{[1]}(3,1)$, $c^{[1]}(2,2)$, $c^{[1]}(3,2)$, $%
c^{[1]}(4,2)$, $c^{[1]}(3,3)$, $c^{[1]}(4,3)$, $c^{[1]}(5,3)$, $c^{[2]}(1,1)$%
, $c^{[2]}(2,2)$, $c^{[2]}(3,3)$ \textit{are positive,}

\smallskip

(ii) \textit{one of following conditions holds; }

\smallskip

\ \ \ (ii-a) $\widehat{c}^{[1]}(3,3)\geq 0$, $c^{[1]}(4,4)\geq 0$ \emph{or},

\smallskip

\ \ \ (ii-b) $\widehat{c}^{[1]}(3,3)\leq 0$, $v_{6}v_{4}v_{2}v_{0}+K\widehat{%
c}^{[1]}(3,3)\geq 0$,

\smallskip

(iii) \textit{one of following conditions holds;}

\smallskip

\ \ \ (iii-a) $\sum_{k=0}^{1}Q_{k}\widehat{c}^{[1]}(3,3-k)\geq 0$, $%
c^{[1]}(5,4)\geq 0$ \emph{or},

\smallskip

\ \ \ (iii-b) $\sum_{k=0}^{1}Q_{k}\widehat{c}^{[1]}(3,3-k)\leq 0$, $%
v_{6}v_{4}v_{2}v_{0}+\sum_{k=0}^{1}K^{k+1}\widehat{c}^{[1]}(3,3-k)\geq 0$,

\smallskip

(iv) \textit{one of following conditions holds;}

\smallskip

\ \ \ (iv-a) $\sum_{k=0}^{2}Q_{k}\widehat{c}^{[1]}(3,3-k)\geq 0$, $%
c^{[1]}(6,4)\geq 0$ \emph{or},

\smallskip

\ \ \ (iv-b) $\sum_{k=0}^{2}Q_{k}\widehat{c}^{[1]}(3,3-k)\leq 0$, $%
v_{6}v_{4}v_{2}v_{0}+\sum_{k=0}^{2}K^{k+1}\widehat{c}^{[1]}(3,3-k)\geq 0$,

\smallskip

(v) \textit{one of following conditions holds;}

\smallskip

\ \ \ (v-a) $\widehat{c}^{[2]}(3,3)\geq 0$, $c^{[2]}(4,4)\geq 0$ \emph{or},

\smallskip

\ \ \ (v-b) $\widehat{c}^{[2]}(3,3)\leq 0$, $v_{7}v_{5}v_{3}v_{1}+K\widehat{c%
}^{[2]}(3,3)\geq 0$.

\medskip

According to Proposition 5.1, we discuss some examples related to the main
theorem (cf. Theorem 4.11).

\medskip

\textbf{Example 5.2.} Let $\alpha $ : $1,1,\sqrt{\frac{106}{100}},\sqrt{x}%
,\left( \sqrt{\frac{111}{100}},\sqrt{\frac{112}{100}},\sqrt{\frac{113}{100}}%
\right) ^{\wedge }$ with $\frac{106}{100}<x<\frac{111}{100}$. Then

\noindent {\footnotesize
\begin{align*}
c^{[1]}(2,1)>0\Longleftrightarrow & -1137824+2385000x-1225625x^{2}>0, \\
c^{[1]}(3,1)>0\Longleftrightarrow & -789728+2068600x-1225625x^{2}>0, \\
c^{[1]}(4,2)>0\Longleftrightarrow & \ 72055895857008-74351127838375x+8426316498750x^{2}>0, \\
c^{[2]}(1,1)>0\Longleftrightarrow & \ {\tiny \frac{50}{47}}<x, \\
c^{[2]}(2,2)>0\Longleftrightarrow & -450800+788127x-342990x^{2}>0, \\
c^{[2]}(3,3)>0\Longleftrightarrow & -1433513345600+2384120972114x \\
& -839363043975x^{2}-126847775250x^{3}>0.
\end{align*}}
Hence condition (i) of Proposition 5.1 holds if and only if $\frac{%
262709-\sqrt{296066681}}{228660}<x<\frac{4(10343+\sqrt{10186611})}{49025}$.

\smallskip

\noindent Since $\widehat{c}^{[1]}(3,3)$ and $c^{[1]}(4,4)$ in (ii-a) of Proposition
5.1 are nonnegative for any $1<x<\frac{111}{100}$ (as is easy to determine because the relevant expression is linear (respectively, quadratic) in $x$), condition (ii) of
Proposition 5.1 holds and since $c^{[1]}(5,4)$ and $v_{6}v_{4}v_{2}v_{0}+%
\sum_{k=0}^{1}K^{k+1}\widehat{c}^{[1]}(3,3-k)$ are nonnegative for any $x\in
(\frac{106}{100},\frac{111}{100})$ in our current example,  condition (iii) of Proposition 5.1
holds automatically.

\noindent Observe that $\sum_{k=0}^{2}Q_{k}\widehat{c}^{[1]}(3,3-k),$ $c^{[1]}(6,4)$,
and $v_{6}v_{4}v_{2}v_{0}+\sum_{k=0}^{2}K^{k+1}\widehat{c}^{[1]}(3,3-k)$ are
nonnegative if and only if certain polynomials $\rho
_{i}(x)=\sum_{j=0}^{3}a_{ij}x^{j}\geq 0$, $i=1,2,3$. One computes that condition (iv) of Proposition 5.1 holds if and only if $%
\frac{106}{100}<x\leq \delta _{1},$ where $\delta _{1}$ is one root of $\rho
_{3}(x)$ on $\frac{106}{100}<x<\frac{111}{100}$ and $\rho
_{3}(x)$ has coefficients {\footnotesize
\begin{align*}
a_{10}& =5009763244345301025528,\ \ a_{11} =17314060626671803639475,\\
a_{12}& =-31888751427147091995000,\ \ a_{13} =10945821827125518750000,\\
a_{20}& =-11515715228709925376854972800,\ \ a_{21} =58140818081665110669480014972,\\
a_{22}& =66850284929099256986049235575,\ \ a_{23} =21452878525521000469568062500,\\
a_{30}& =-287892880717748134421374320000+15930706453117441780709304096\sqrt{7}, \\
a_{31}& =1453520452041627766737000374300+55057535436693721890883015700\sqrt{7}, \\
a_{32}& =-1671257123227481424651230889375-101404061103230706541844340000\sqrt{7}, \\
a_{33}& =536321963138025011739201562500+34806969094374905089725000000\sqrt{7}.
\end{align*}}
Observe that $\widehat{c}^{[2]}(3,3)$, $c^{[2]}(4,4)$, and $%
v_{7}v_{5}v_{3}v_{1}+K\widehat{c}^{[2]}(3,3)$ are nonnegative if and only if
certain polynomials $\rho _{i}(x)=\sum_{j=0}^{3}a_{ij}x^{j}\geq 0$, $i=4,5,6$, for
some $a_{ij}\in \mathbb{R}$. One shows that condition (v) of Proposition 5.1 holds
if and only if $\frac{106}{100}<x\leq \delta _{2}$, where $\delta _{2}$ is
one root of $\rho _{6}(x)$ on $\frac{106}{100}<x<\frac{111}{100}$ and $\rho _{6}(x)$ has coefficients {\footnotesize
\begin{align*}
a_{40}& = -901600,\ \ a_{41} = 1582879,\ \ a_{42} = -693750,\ \ a_{43} = 0,\\
a_{50}& = -32970806948800,\ \ a_{51} = 54940864614622,\\
a_{52}& = - 19516290706245,\ \ a_{53} = -  2816020610550,\\
a_{60}& = -161987008052800-11468106764800\sqrt{7}, \\
a_{61}& = 269670875488882+20133790336912\sqrt{7}, \\
a_{62}& = -95375375706225-8824311300000\sqrt{7}, \\
a_{63}& = -14080103052750.
\end{align*}}
Hence, combining the conditions above , we get that $W_{\alpha }$ is
semi-cubically hyponormal with p.d.c. if and only if $\frac{262709-\sqrt{%
296066681}}{228660}<x<\frac{4(10343+\sqrt{10186611})}{49025}$.

\medskip

\textbf{Example 5.3.} Let $\alpha $ : $1,1,\sqrt{x},\sqrt{\frac{%
109}{100}},\left( \sqrt{\frac{111}{100}},\sqrt{\frac{112}{100}},\sqrt{\frac{%
113}{100}}\right) ^{\wedge }$ with $1<x<\frac{109}{100}$. Then

\noindent {\footnotesize
\begin{align*}
c^{[1]}(2,1)>0 \Longleftrightarrow & -1404816 + 2485825 x - 1090000 x^2>0,\\
c^{[1]}(3,1)>0 \Longleftrightarrow & -1456336 + 2537345 x - 1090000 x^2>0,\\
c^{[1]}(3,2)>0 \Longleftrightarrow & -10994090016 + 19612265981 x -8666453750 x^2>0,\\
c^{[1]}(4,2)>0 \Longleftrightarrow & -140086383317216 + 246270929381981 x -106779376653750 x^2>0,\\
c^{[1]}(4,3)>0 \Longleftrightarrow & -987632956940439264 + 1763044935312612532 x\\
                                   & -772409239899810625 x^2 -6347196840234375 x^3 > 0,\\
c^{[1]}(5,3)>0 \Longleftrightarrow & -12598978259232068585504 + 22152829244755615420532 x\\
                                   & -9516854244805566710625 x^2 - 78203812268527734375 x^3 > 0,\\
c^{[2]}(1,1)>0 \Longleftrightarrow & 1<x<\frac{118}{109}.
\end{align*}}
\noindent Note that $c^{[2]}(2,2), c^{[2]}(3,3)$ are positive, for $1<x<\frac{109}{100}$. Hence condition (i) of Proposition 5.1 holds if and only if $\frac{99433 - \sqrt{86925073}}{87200} < x< \frac{118}{109}$.

\smallskip

\noindent Since $\widehat{c}^{[1]}(3,3)$ and $c^{[1]}(4,4)$ in (ii-a) of Proposition
5.1 are nonnegative for any $1<x<\frac{111}{100}$ (as is easy to determine because the relevant expression is linear (respectively, quadratic) in $x$), condition (ii) of
Proposition 5.1 holds.

\noindent Observe that $\sum_{k=0}^{1}Q_{k}\widehat{c}^{[1]}(3,3-k),$ $c^{[1]}(5,4)$,
and $v_{6}v_{4}v_{2}v_{0}+\sum_{k=0}^{1}K^{k+1}\widehat{c}^{[1]}(3,3-k)$ are
nonnegative if and only if certain polynomials $\tau_{i}(x)=\sum_{j=0}^{3}b_{ij}x^j \geq 0$, $i=1,2,3$. One computes that condition (iii) of Proposition 5.1 holds if and only if $\delta_3%
\leq x < \frac{109}{100},$ where $\delta_3$ is one root of $\tau_{3}(x)$ on $1<x<\frac{109}{100}$ and $\tau
_{3}(x)$ has coefficients {\footnotesize
\begin{align*}
b_{10}& = -141090422420062752,\ \ b_{11} = 250594246445105815,\ \ b_{12} = -110492596718511250,\\
b_{13}& = 0,\ \ b_{20} = -6559637294111604097794672,\ \ b_{21} = 11705749677041669653934138,\\
b_{22}& = - 5124944374775526264385625,\ \ b_{23} = -  43715931058107003515625,\\
b_{30}& = -6559637294111604097794672 - 897315898294149974185728 \sqrt{7}, \\
b_{31}& = 11705749677041669653934138 +1593745326573356449009160 \sqrt{7}, \\
b_{32}& = -5124944374775526264385625 - 702717888136577832470000 \sqrt{7}, \\
b_{33}& =- 43715931058107003515625.
\end{align*}}

\noindent Observe that $\sum_{k=0}^{2}Q_{k}\widehat{c}^{[1]}(3,3-k),$ $c^{[1]}(6,4)$,
and $v_{6}v_{4}v_{2}v_{0}+\sum_{k=0}^{2}K^{k+1}\widehat{c}^{[1]}(3,3-k)$ are
nonnegative if and only if certain polynomials $\tau
_{i}(x)=\sum_{j=0}^{3}b_{ij}x^{j}\geq 0$, $i=4,5,6$. One computes that condition (iv) of Proposition 5.1 holds if and only if $\delta_4%
\leq x < \frac{109}{100} ,$ where $\delta_4$ is one root of $\tau_{6}(x)$ on $1<x<\frac{109}{100}$ and $\tau
_{6}(x)$ has coefficients {\footnotesize
\begin{align*}
b_{40}& = -1828509633963457884992,\ \ b_{41} = 3177706249776013464215,\\
b_{42}& = -1361379284168777111250,\ \ b_{43} = 0,\\
b_{50}& = -83769913286973551581368596752,\ \ b_{51} = 147175163957054889298562957338,\\
b_{52}& = - 63144439641609259103495285625,\ \ b_{53} = -  538623986566936390316015625,\\
b_{60}& = -83769913286973551581368596752 - 11629072594697373118276761088 \sqrt{7}, \\
b_{61}& = 147175163957054889298562957338 +20209779580525476094576266760 \sqrt{7}, \\
b_{62}& = -63144439641609259103495285625 - 8658187099730775473862870000 \sqrt{7}, \\
b_{63}& = - 538623986566936390316015625.
\end{align*}}

\noindent Since $c^{[2]}(4,4)$ and $v_{7}v_{5}v_{3}v_{1}+%
K\widehat{c}^{[2]}(3,3)$ are nonnegative for any $x\in
(1,\frac{109}{100})$ in our current example,  condition (v) of Proposition 5.1
holds automatically.

\noindent Hence, combining above conditions, we get that $W_{\alpha }$ is semi-cubically
hyponormal with p.d.c. if and only if $\frac{99433 - \sqrt{86925073}}{87200} < x< \frac{118}{109}$.

\medskip

We give next an example with two variables$\ x$ and $y$ which generalizes
Examples 5.2 and 5.3.

\medskip

\textbf{Example 5.4.} Let $\alpha :=\alpha (x,y):1,1,\sqrt{x},\sqrt{y}%
,\left( \sqrt{\frac{111}{100}},\sqrt{\frac{112}{100}},\sqrt{\frac{113}{100}}%
\right) ^{\wedge }$ with $1<x<y<\frac{111}{100}$. According to Proposition
3.21, $W_{\alpha }$ is semi-cubically hyponormal with p.d.c. if and only if
the following conditions hold:

(a) the coefficients $c^{[1]}(2,1)$, $c^{[1]}(3,1),$ $c^{[1]}(3,2),$ $%
c^{[1]}(4,2),$ $c^{[1]}(4,3),$ $c^{[1]}(5,3),$ $c^{[2]}(1,1),$ $c^{[2]}(2,2),
$ $c^{[2]}(3,3)$ are positive, equivalently, polynomials $p_{j}(x,y)$ $%
(1\leq j\leq 9)$ with degrees $3$ or less in $x$ [and $y$] (which are in
Appendix) are positive.

Note that since $\widehat{c}^{[1]}(3,3)$ and $c^{[1]}(4,4)$ in (ii-a) of Proposition
3.21 are nonnegative for any $1<x<y<\frac{111}{100},$ a condition (ii) of Proposition
3.21 holds.

(b) either $\phi _{1}(x,y)\geq 0$ and $\phi _{2}(x,y)\geq 0$ for $1<x<y<%
\frac{111}{100}$ or $\phi _{1}(x,y)\leq 0$ and $\phi _{3}(x,y)\geq 0$ for $%
1<x<y<\frac{111}{100}.$ (Note that $\phi _{1}(x,y)\geq 0$, $\phi
_{2}(x,y)\geq 0$ and $\phi _{3}(x,y)\geq 0$ are equivalent to $%
\sum_{k=0}^{1}Q_{k}\widehat{c}^{[1]}(3,3-k)\geq 0$, $c^{[1]}(5,4)\geq 0$ and
$v_{6}v_{4}v_{2}v_{0}+\sum_{k=0}^{1}K^{k+1}\widehat{c}^{[1]}(3,3-k)\geq 0,$
respectively, in (iii) of Proposition 3.21.)

(c) either $\phi _{4}(x,y)\geq 0$ and $\phi _{5}(x,y)\geq 0$ for $1<x<y<%
\frac{111}{100}$ or $\phi _{4}(x,y)\leq 0$ and $\phi _{6}(x,y)\geq 0$ for $%
1<x<y<\frac{111}{100}.$ (Note that $\phi _{4}(x,y)\geq 0$, $\phi
_{5}(x,y)\geq 0$ and $\phi _{6}(x,y)\geq 0$ are equivalent to $%
\sum_{k=0}^{2}Q_{k}\widehat{c}^{[1]}(3,3-k)\geq 0$, $c^{[1]}(6,4)\geq 0$ and
$v_{6}v_{4}v_{2}v_{0}+\sum_{k=0}^{2}K^{k+1}\widehat{c}^{[1]}(3,3-k)\geq 0,$
respectively, in (iv) of Proposition 3.21.)

(d) either $\phi _{7}(x,y)\geq 0$ and $\phi _{8}(x,y)\geq 0$ for $1<x<y<%
\frac{111}{100}$ or $\phi _{7}(x,y)\leq 0$ and $\phi _{9}(x,y)\geq 0$ for $%
1<x<y<\frac{111}{100}.$ (Note that $\phi _{7}(x,y)\geq 0$, $\phi
_{8}(x,y)\geq 0$ and $\phi _{9}(x,y)\geq 0$ are equivalent to $%
\widehat{c}^{[2]}(3,3)\geq 0$, $c^{[2]}(4,4)\geq 0$ and
$v_{7}v_{5}v_{3}v_{1}+K\widehat{c}^{[2]}(3,3)\geq 0,$
respectively, in (v) of Proposition 3.21.)

\medskip

In the diagram that follows positivity of the coefficients $c^{[1]}(2,1)$, $c^{[1]}(3,1),$ $c^{[1]}(3,2),$ $%
c^{[1]}(4,2),$ $c^{[1]}(4,3),$ $c^{[1]}(5,3),$ $c^{[2]}(1,1),$ $c^{[2]}(2,2),
$ $c^{[2]}(3,3)$ in Example 5.4(i) can be represented by regions bounded by polynomials $%
p_{i}(x,y),$ $1\leq i\leq 9.$ Figure 2 shows the behaviors of the relevant polynomials and the relevant $\phi_i$ for
this example and the region of the semi-cubic hyponormality with p.d.c.  (Functions not bounding the final region are not displayed.)

\begin{center}
\includegraphics[width=10cm]{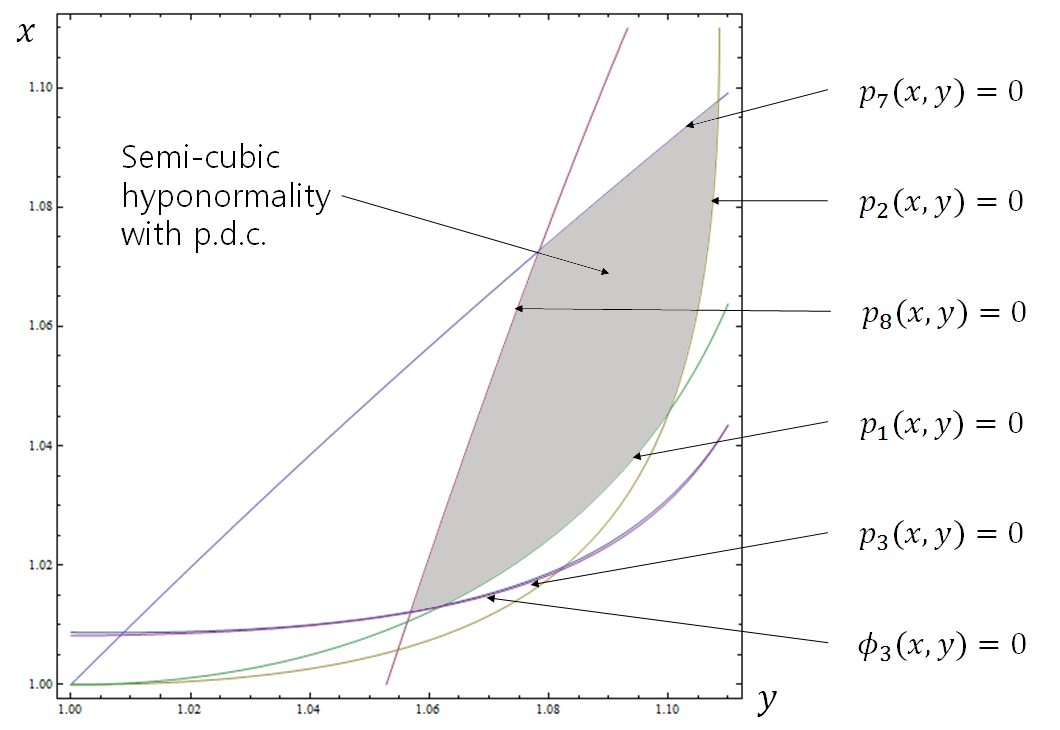}\\[0pt]
\text{Figure 2}
\end{center}

Continuing this same example, we provide for comparison the parameter region in which the weighted shift is quadratically hyponormal with positive determinant coefficients (i.e., positively quadratically hyponormal), using the techniques of Theorem 3.6 of \cite{EJLP}.  Conveniently, Figure 1 of that paper happens to indicate the coefficients which must be checked for a backward $4$-step extension.  In the following diagram, we have been content to approximate the region with only three of the thirty-seven relevant bounding curves, and, in particular, have for clarity omitted a few bounding curves which make very minor changes to the region.  It is clear that there is (substantial) overlap of the regions of quadratic hyponormality with p.d.c. and semi-cubic hyponormality with p.d.c., and $y = \frac{108}{100}$, $x = \frac{105}{100}$ is a (sample) point in common. These parameter values provide an example of a shift with both of these properties, but which is not cubically hyponormal, since cubically hyponormal weighted shifts have the flatness property in which the first two weights equal forces all weights equal (see \cite{LCL}).

\medskip

\begin{center}
\includegraphics[width=10cm]{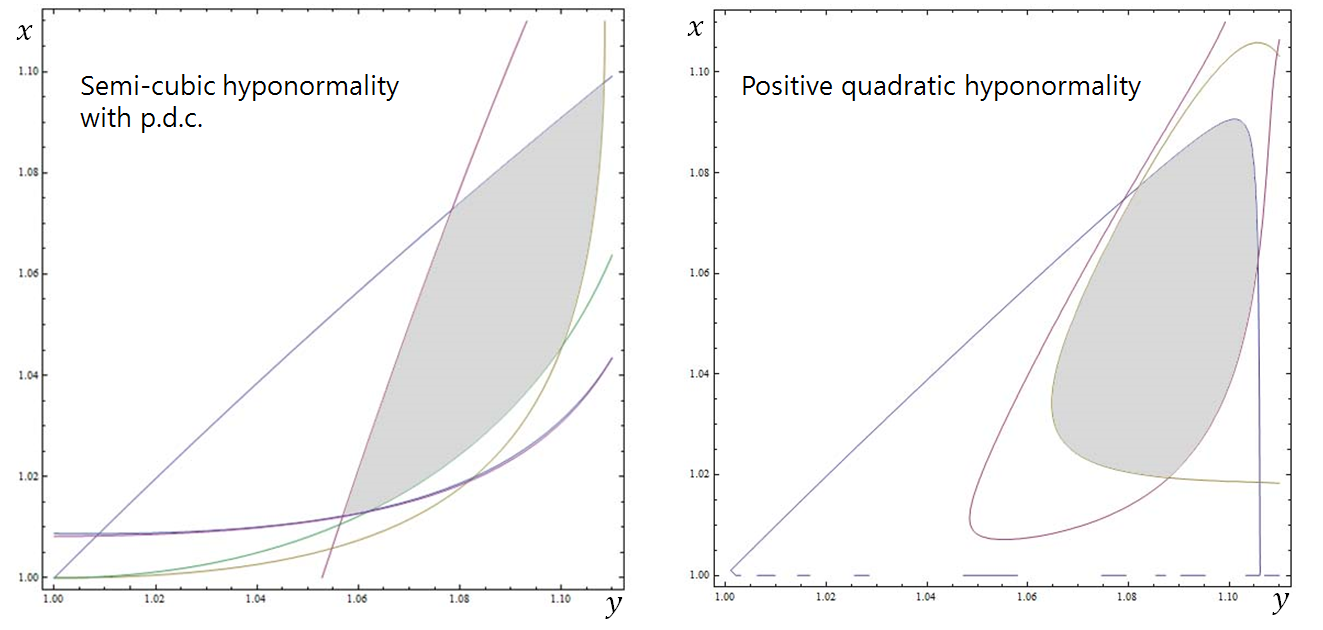}\\[0pt]
\text{Figure 3}
\end{center}

\medskip

We provide for the reader's convenience both a side-by-side and a single picture of the two regions;  it is clear that neither is contained in the other.  (This assertion is unaffected by our slightly simplified version of the positive quadratic hyponormality region.)

\begin{center}
\includegraphics[width=8cm]{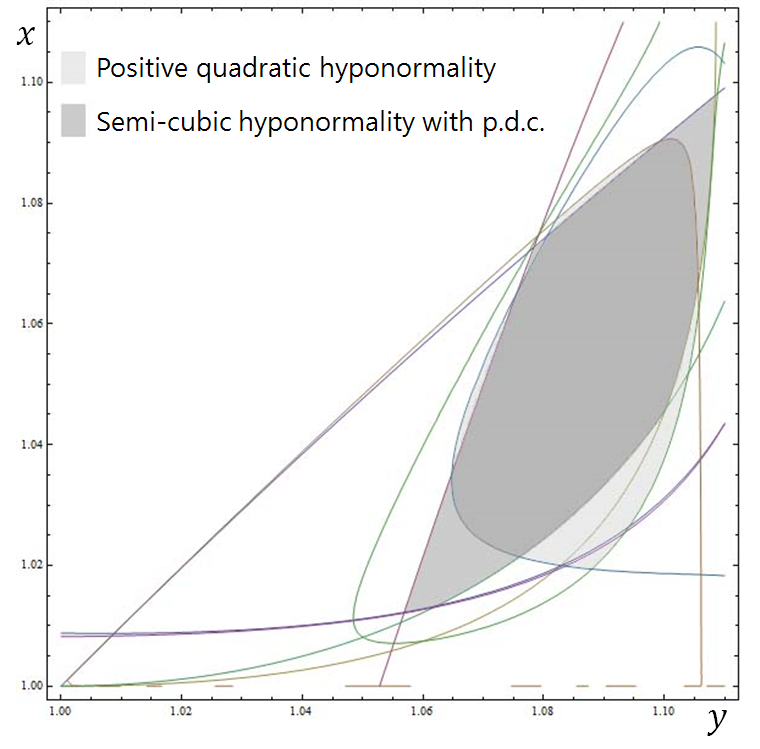}\\[0pt]
\text{Figure 4}
\end{center}

\medskip

Observe also that $\alpha (x,y)$ above is a backward $3$-step extension with Stampfli tail if and only if $y=%
\frac{3108}{2825}$;   that is, the weight sequence must be of the form $1,1,\sqrt{x},\sqrt{\frac{3108}{2825}},\left( \sqrt{%
\frac{111}{100}},\sqrt{\frac{112}{100}},\sqrt{\frac{113}{100}}\right)
^{\wedge }=1,1,\sqrt{x}, \left(\sqrt{\frac{3108}{2825}}, \sqrt{\frac{111}{100%
}},\sqrt{\frac{112}{100}}\right) ^{\wedge }$.

Note that, in \cite[Example 4.1]{BEJL}, one considered a backward $3$-step
extension weight sequence $\alpha :1,1,\sqrt{x},\left( \sqrt{\frac{111}{100}}%
,\sqrt{\frac{112}{100}},\sqrt{\frac{113}{100}}\right) ^{\wedge }$ with $1<x<%
\frac{111}{100}$ and discussed the equivalent conditions of semi-cubic
hyponormality with p.d.c. of $W_{\alpha }$. The weight sequence $\alpha :1,1,%
\sqrt{x},\left( \sqrt{\frac{111}{100}},\sqrt{\frac{112}{100}},\sqrt{\frac{113%
}{100}}\right) ^{\wedge }$ can be put in the form of a backward $4$-step
extension weight sequence as $\alpha :1,1,\sqrt{x},\sqrt{\frac{111}{100}}$, $%
\left( \sqrt{\frac{112}{100}},\sqrt{\frac{113}{100}},\sqrt{\frac{644}{565}}%
\right) ^{\wedge }$.

\bigskip

\textbf{Appendix: The expressions of polynomials in Example 5.4.}

\smallskip

\noindent {\small$p_1(x,y):=-87801 + 18426 x + 201250 x y - 62500 x^2 y - 69375 x y^2$,}

\noindent {\small$p_2(x,y):=-435897 + 366522 x + 316400 y - 115150 x y - 62500 x^2 y - 69375 x y^2$,}

\noindent {\small$p_3(x,y):=-8123260719 + 7850825094 x + 6822137700 y - 609514150 x y$}

\quad \quad {\small$  -5695718750 x^2 y- 5016991875 x y^2 + 4769531250 x^2 y^2$,}

\noindent {\small$p_4(x,y):=-131357230786383 + 128000551450758 x + 112478744797300 y$}

\quad \quad {\small$ - 35933010037750 x y - 70176950718750 x^2 y - 61814356891875 x y^2$}

\quad \quad {\small$ + 58765394531250 x^2 y^2$,}

\noindent {\small$p_5(x,y):=-729737521907523801 + 647639252666034426 x + 612853634952978300 y$}

{\small$+ 61743428759974150 x y - 473135879119093750 x^2 y - 509005371229162500 x y^2$}

{\small$+ 348191411957812500 x^2 y^2 - 333894287109375 x^3 y^2 + 41509737773437500 x^2 y^3$,}

\noindent {\small$p_6(x,y):=-11898782682983018608569 + 10887249907658628019194 x$}

\quad \quad {\small$+ 10193895909890838827500 y - 1882185487883551691050 x y$}

\quad \quad {\small$ - 5829507166626354093750 x^2 y - 6271455178914511162500 x y^2$}

\quad \quad {\small$ + 4290066386732207812500 x^2 y^2 - 4113911511474609375 x^3 y^2$}

\quad \quad {\small$ + 511441479106523437500 x^2 y^3$,}

\noindent {\small$p_7(x,y):=-1 + 2 y - x y$,}

\noindent {\small$p_8(x,y):=-90160 + 180320 y - 21410 x y - 86247 y^2 + 16650 x y^2$,}

\noindent {\small$p_9(x,y):=-71675667280 + 143351334560 y - 22778571655 x y - 61406571072 y^2$}

\quad \quad {\small$+ 13236466950 x y^2 + 4812890625 x^2 y^2 - 5983385625 x y^3$,}

\noindent {\small$\phi_1(x,y):=-833985739322884344+775904394232181844x+700404154231975200y$}

\quad \quad {\small$+866850580839425xy-567462345089140625x^{2}y-548398399244797500xy^{2}$}

\quad \quad {\small$+448253315343750000x^{2}y^{2}+23719850156250000x^2y^{3}$,}

\noindent {\small$\phi_{2}(x,y):=-9693506945021586690580521+8631115634040960029101146x$}

\quad \quad {\small$+8140873654364803833354300y+765254451524447878242550xy$}

\quad \quad {\small$-6235415266776031685781250x^{2}y-6735143739224307305482500xy^{2}$}

\quad \quad {\small$+4600718666945449248046875x^{2}y^{2}-54373568447159912109375x^{3}y^{2}$}

\quad \quad {\small$+532706793238642148437500x^{2}y^{3}+45664417777368164062500x^{3}y^{3}$,}

\noindent {\small$\phi _{3}(x,y):=-9693506945021586690580521-1326008970008249128892304\sqrt{7}$}

\quad \quad {\small$+(8631115634040960029101146+1233661606079765237777304\sqrt{7})x$}

\quad \quad {\small$+(8140873654364803833354300+1113618791487596680843200\sqrt{7})y$}

\quad \quad {\small$+(765254451524447878242550x+1378262950614937209550\sqrt{7}xy)$}

\quad \quad {\small$-(6235415266776031685781250+902245834972000562968750\sqrt{7})x^{2}y$}

\quad \quad {\small$-(6735143739224307305482500+871934809253653701885000\sqrt{7})xy^{2}$}

\quad \quad {\small$+(4600718666945449248046875+712707530783840812500000\sqrt{7})x^{2}y^{2}$}

\quad \quad {\small$+(532706793238642148437500+37713755273532187500000\sqrt{7})x^{2}y^{3}$}

\quad \quad {\small$-54373568447159912109375x^{3}y^{2}+45664417777368164062500x^{3}y^{3}$,}

\noindent {\small$\phi _{4}(x,y):=-15147506089410954390264 + 14431885836548408887764 x$}

\quad \quad {\small$ +13058028690302041695200 y - 4417668640003352700575 x y$}

\quad \quad {\small$ -6991703553843301640625 x^2 y - 6756816677095149997500 x y^2$}

\quad \quad {\small$ +5522929098350343750000 x^2 y^2 + 292252273775156250000 x^2 y^3$,}

\noindent {\small$\phi _{5}(x,y):=-6371142354921078097021338225 + 5847553421217386053177843050 x$}

\quad {\small$ +5460819072005697558820637836 y - 1071522896299238745277246282 x y$}

\quad {\small$ -3073062060077899456020431250 x^2 y - 3319348240439307612433995300 x y^2$}

\quad {\small$ +2267418187817395207407421875 x^2 y^2 - 26797469473498291083984375 x^3 y^2$}

\quad {\small$ +262539215979732396435937500 x^2 y^3 + 22505251657398125976562500 x^3 y^3$,}

\noindent {\small$\phi _{6}(x,y):=-159278558873026952425533455625-24084019666956377508070491024\sqrt{7}$}

{\small$+(146188835530434651329446076250+22946207795993527485642576024\sqrt{7})x$}

{\small$+(136520476800142438970515945900+20761821644604776025950363200\sqrt{7})y$}

{\small$-(26788072407480968631931157050+7023942936871570679922430450\sqrt{7})xy$}

{\small$-(76826551501947486400510781250+11116570932690018936337968750\sqrt{7})x^{2}y$}

{\small$-(82983706010982690310849882500+10743108784814267260925085000\sqrt{7})xy^{2}$}

{\small$+(56685454695434880185185546875+8781269486787702650812500000\sqrt{7})x^{2}y^{2}$}

{\small$+(6563480399493309910898437500+464671178725190082187500000\sqrt{7})x^{2}y^{3}$}

{\small$-669936736837457277099609375x^{3}y^{2}+562631291434953149414062500x^{3}y^{3}$,}

\noindent {\small$\phi _{7}(x,y):=-90160 + 180320 y - 20785 x y - 87024 y^2 + 16650 x y^2$,}

\noindent {\small$\phi _{8}(x,y):=-1648540347440 + 3297080694880 y - 518903268065 x y - 1418571958272 y^2$}

\quad \quad {\small$ +304438739850 x y^2 + 106846171875 x^2 y^2 - 132831160875 x y^3$,}

\noindent {\small$\phi _{9}(x,y):=-8099350402640-573405338240\sqrt{7}+(16198700805280+1146810676480\sqrt{7})y$}

{\small$-(2561468897015+132189773240\sqrt{7})xy-(6954494590176+553460804736\sqrt{7})y^{2}$}

{\small$+(1495720765350+105891735600\sqrt{7})xy^{2}+534230859375x^{2}y^{2}-664155804375xy^{3}$.}

\bigskip

\bigskip

\noindent Seunghwan Baek

\noindent Kyungpook National University, Daegu 702-701, Korea

\noindent E-mail: seunghwan@knu.ac.kr

\medskip

\noindent George R. Exner

\noindent Bucknell University, Lewisburg, Pennsylvania 17837, USA

\noindent E-mail: exner@bucknell.edu

\medskip

\noindent Il Bong Jung

\noindent Kyungpook National University, Daegu 702-701, Korea

\noindent E-mail: ibjung@knu.ac.kr

\medskip

\noindent Chunji Li

\noindent Northeastern University, Shenyang 110004, P. R. China

\noindent E-mail: lichunji@mail.neu.edu.cn

\end{document}